\documentclass[a4paper]{article}
\newcommand{\rank}{\hbox{rank}\ }
\newcommand{\Stab}{\hbox{Stab}\ }
\newcommand{\g}{\pi (g)}
\newcommand{\h}{\pi (h)}
\newcommand{\B}{\cal O}
\usepackage{amsmath}
\usepackage{amssymb}
\usepackage{latexsym}

\newtheorem{Theorem}{Theorem}
\newtheorem{Proposition}{Proposition}
\newtheorem{Lemma}{Lemma}
\newtheorem{Corollary}{Corollary}
\newtheorem{Remark}{Remark}

\begin{document}

\title{The structure of 3-manifolds with 2-generated fundamental
group}

\author{Michel Boileau and Richard Weidmann} \date{} \maketitle

The main purpose of this article is to describe compact orientable
irreducible
3-manifolds that have a non-trivial JSJ-decomposition and 2-generated
fundamental
group. The rank of a group is the minimal number of elements needed
to generate it.
A natural question is whether the Heegaard genus of such a manifold
is equal
to $2$.

When the JSJ-decomposition is empty, there are examples of closed
3-manifolds of Heegaard
genus $3$ that have 2-generated fundamental group \cite{BZg}. These
are
Seifert fibered manifolds. Furthermore the second named author
\cite{W2} has recently
found graph manifolds with the same property. At this point these are
the only known
examples of 3-manifolds that have 2-generated fundamental group but
are not of Heegaard genus $2$. In particular, it is still an open
problem to
find such examples that admit a complete hyperbolic structure. There
are also examples
of Seifert manifolds with Heegaard genus $g+1\ge 3$ and fundamental
group of rank $g$ \cite{MSc}.

When the JSJ-decomposition is non-trivial, our main result with
respect to this question
is the following; we will denote the base spaces by their topological
type, followed
by a list with the orders of their cone points. We denote the
M\"obius band by $M\ddot o$, the disk by $D$, the annulus by $A$ and
the 2-punctured
disk by $\Sigma$. We furthermore denote by $Q^3$ the orientable
circle bundle over
the M\"obius band.

\begin{Theorem}\label{heegrank} Let $M$ be a compact, orientable,
irreducible 3-manifold with\break $\rank \pi_1(M)=2$. If $M$
has a non-trivial JSJ-decomposition then one of the following holds:

\begin{enumerate} \item $M$ is of Heegaard genus $2$.

\item $M=S\cup _T H$ where $S$ is a Seifert manifold with basis
$D(p,q)$ or $A(p)$, $H$ is a hyperbolic manifold and $\pi_1(H)$ is
generated by a pair of elements with a single parabolic element. The
gluing map identifies the fibre of $S$ with the curve corresponding
to the parabolic generator of $\pi_1(H)$.

\item $M=S_1\cup _T S_2$ where $S_1$ is a Seifert manifold over
$M\ddot o$ or $M\ddot o(p)$ and $S_2$ is a Seifert manifold over
$D(2,2l+1)$. The gluing map identifies the fibre of $S_1$ with a
curve on the boundary of $S_2$ that has intersection number one with
the fibre of $S_2$.

\item $M= Q^3\cup H$ where $H$ is a hyperbolic manifold that admits a
finite-sheeted irregular covering by the exterior of a hyperbolic
2-bridge link.

\end{enumerate} \end{Theorem}

\begin{Remark} It follows from T. Kobayashi's work \cite{Ko1} that
the three manifolds
of type 3) do not Heegaard genus two unless either $S_1$ is the exterior of
a 1-bridge knot in a lens space and the meridian is glued with the fibre of
$S_2$,
or $S_2$ is the exterior of a 2-bridge knot in $\mathbb S^3$ and the
meridian is
glued with the fibre of $S_1$ (see \cite{W2}). We do not yet
know an example of
manifolds of type 2) that does not have Heegaard genus two. Moreover
we conjecture that
there is no example of type 4).
\end{Remark}

It is known that a Heegaard genus 2 closed 3-manifold is a 2-fold
branched covering of the 3-sphere $\mathbb S^{3}$ \cite{BH}. For a
closed 3-manifolds with
a 2-generated fundamental group, that is geometric or has a
non-trivial JSJ-decomposition, we obtain the following related
result:

\begin{Corollary}\label{branchcov} Let $M^3$ be a closed, orientable,
irreducible
3-manifold with \break $\rank \pi_1(M^3)=2$. If $M$ is geometric or
has a non-trivial JSJ-decomposition, then
$M$ is a 2-fold branched covering of a homotopy sphere.
\end{Corollary}

\begin{Remark} In general, we do not know how to show that
these homotopy
spheres are the true sphere $\mathbb S^{3}$.
\end{Remark}

With respect to Thurston's geometrization conjecture, we have the
following corollary for closed 3-manifolds containing a two-generator
knot:

\begin{Corollary}\label{geom} Let $M^3$ be a closed, orientable,
irreducible 3-manifold. If there is a two-generator knot $k\subset M$
that is not in a ball, then either $M$ is geometric or has a
non-trivial JSJ-decomposition. \end{Corollary}

It follows from Thurston's orbifold theorem ([BP],[CHK]) that a Heegaard
genus two
closed 3-manifolds has a geometric decomposition in the sens of Thurston.
Since by [JS1, Lemma 5.4], a compact, orientable, irreducible 3-manifold with
rank two fundamental group is either a handlebody of genus two or has
Euler characteristic zero,
Corollary \ref{geom} gives the following algebraic version of this result:

\begin{Corollary}\label{homgeom} Let $M^3$ be a closed, orientable,
irreducible 3-manifold
which is the union of two non-empty,
compact, connected, irreducible, orientable  3-manifolds $U$ and $V$ along
their boundaries.
If both $1\leq \rank \pi_1(U) \leq 2$ and $1\leq \rank \pi_1(V) \leq 2$,
then either $M^3$ is geometric or has a
non-trivial JSJ-decomposition. \end{Corollary}

  F.H. Norwood \cite{N} has shown that a two-generator knot in
$\mathbb S^{3}$
  is prime. Later on M. Scharlemann \cite{Scha}
has shown that a
tunnel number one
(i.e. Heegaard genus two) knot
in $\mathbb S^{3}$ is Conway irreducible (i.e doubly prime). The
following corollary extends this last result to the case of a
two-generator
knot in $\mathbb S^{3}$.

\begin{Corollary}\label{conway} A two-generator knot $k\subset
\mathbb
S^{3}$ is prime and Conway irreducible. \end{Corollary}

\noindent Theorem~\ref{heegrank} follows from the classification of
3-manifolds with non-trivial JSJ-decomposition and Heegaard genus 2
as given by T. Kobayashi's work \cite{Ko1}, \cite{Ko2}, \cite{Ko3}
and the following theorem that describes (almost) precisely the
JSJ-decomposition of a compact orientable 3-manifold with 2-generated
fundamental group:

\begin{Theorem}\label{JSJ} Let $M^3$ be a compact, orientable,
irreducible 3-manifold with non-trivial JSJ-decomposition. If
$\rank(\pi_1(M^3))=2$ then the
JSJ-decomposition is of one of the following types. If $M^3$ is of
one of the types (1)-(9) the converse holds:

\begin{enumerate}

\item $M_1$ is a Seifert manifold with base $A(p)$ or $D(p,q)$ and
$M_2$ is a hyperbolic manifold whose fundamental group has a
generating pair $\{g,h\}$ where
$g$ is parabolic and $M^3$ is obtained from $M_1$ and $M_2$ by gluing
boundary components such that the fibre of $M_1$ gets identified with
the curve corresponding to the parabolic generator.

\item $M_1$ is a Seifert manifold with base $A(p)$ or $D(p,q)$ and
$M_2$ is the exterior of a 1-bridge knot in a lens space which is
Seifert and where the meridian is not the fibre. $M^3$ is obtained
from $M_1$ and $M_3$ by gluing boundary components such that the
fibre of $M_1$ gets identified with the meridian curve of $M_2$.

\item $M_1$ is a Seifert manifold with base $A(p)$ or $D(p,q)$ and
$M_2$ is a Seifert space with base space $A(p)$. $M^3$ is obtained
from $M_1$ and $M_3$ by gluing boundary components such that the
fibre of $M_1$ has intersection number 1
with the fibre of $M_2$.

\item $M_1$ is a Seifert space with base space $M\ddot o$, $M\ddot
o(p)$, $M\ddot
o(p,q)$, $D(p,q)$, $D(p,q,r)$, $A(p)$, $A(p,q)$, $\Sigma$,
$\Sigma(p)$ or the once punctured M\"obius band with at most one cone
point. and $M_2$ is the exterior of a 2-bridge knot. Boundary
components are glued such that the fibre of
$M_1$ is identified with the meridian curve of $M_2$.

\item $M_1$ is a Seifert space with base space $M\ddot o$ or $M\ddot
o(p)$ and $M_2$ is a Seifert space over $D(2,2l+1)$. Boundary
components are glued such that the fibre of $M_1$ has intersection
number 1 with the fibre of $M_2$.

\item $M_1$ and $M_2$ are Seifert space with base of type $A(p)$ or
$D(p,q)$ and
$M_3$ is the exterior of a 2-bridge link. Boundary components are
glued such that
the fibres of $M_1$ and $M_2$ get identified with the meridian curves
of $M_3$.

\item $M_1$ is a Seifert space with base space of type $A(p)$,
$\Sigma$, $A(p,q)$
or $\Sigma(p)$ and $M_2$ is the exterior of a 2-bridge link. Two
boundary components of $M_1$ are glued to the boundary components of
$M_2$ such that fibre
of $M_1$ gets identified with the meridian curves of $M_2$.

\item $M$ is the exterior of a 2-bridge link and $M^3$ is obtained
from $M$ by identifying the boundaries such that the meridians get
identified.

\item $M$ is a Seifert space with base $A(p)$, $\Sigma$, $A(p,q)$ or
$\Sigma(p)$
and two boundary are identified such that the fibre in one component
has intersection number $1$ with the fibre in the other.

\item $M$ is obtained from $M_1$ and $M_2$ by gluing along their
boundary where $M_1= Q^3$ and $M_2$ is a hyperbolic manifold that
admits a finite-sheeted irregular covering by the exterior of a
hyperbolic 2-bridge link. \end{enumerate}
\end{Theorem}

For the last case see the discussion in later chapters on the
restrictions of the
gluing map. We conjecture that this case does not occur.

We give some results which are not direct consequences of Theorem
\ref{heegrank}
or Theorem \ref{JSJ} but which are immediate consequences of their
proof.

\begin{Corollary}\label{twogen} Let $M^{3}$ be a compact orientable
$3$-manifold. Any $2$-generated subgroup $U$ of $\pi_{1}(M^3)$ is of
one of the following types:

\begin{enumerate} \item $U$ is of finite index in $\pi_{1}(M^3)$
\item $U$ is free or free Abelian \item $U$ is the fundamental group
of a Seifert fibered manifold \item $U$ is a lattice in $PSL(2,\Bbb
C)$
\item $U$ is the
fundamental group of one of the 3-manifolds described in
Theorem~\ref{JSJ}
\end{enumerate}

If $M^3$ is not closed then (1) implies one of the cases (2)-(5)
\end{Corollary}

\begin{Corollary}\label{twoperiph} Let $M^3$ be a compact,
orientable, irreducible
3-manifold with incompressible boundary. If $\pi_1(M^3)$ is generated
by two peripheral elements
then $M^3$ is homeomorphic to the exterior of a 2-bridge knot or link
in $ \mathbb S^3$. \end{Corollary}

\begin{Corollary}\label{oneperiph} Let $M^3$ be a compact orientable,
irreducible
3-manifold with incompressible boundary. If $\pi_1(M^3)$ is generated
by two elements one of which is peripheral then either $M^3$ has
Heegaard genus $2$ or $M^3$ is hyperbolic. \end{Corollary}

\begin{Remark} A stronger version of Corollary \ref{oneperiph} has been
obtained by
S. Bleiler and A. Jones [BJ2] in the case of a knot exterior in  $\mathbb
S^{3}$
(cf. also [B3]).
\end{Remark}

Tunnel number one satellite knots has been classified by Morimoto and
Sakuma \cite{MSa}:
their exteriors are obtained by gluing a torus knot exterior to
a 2-bridge link exterior along a boundary component, such that the
gluing map
identifies a regular fibre of the torus knot exterior with a meridian
of
the link exterior.

The following corollary shows that the existence of a two generator
satellite knot
which has tunnel number at least two reduces to the existence of a two
generator, but not 2-bridge,
  hyperbolic link, with a single meridional generator and with both
components unknotted.

If such a link has tunnel number one, it has been shown that it is a
2-bridge link. (cf. [Ku]).

\begin{Corollary}\label{sat} Let $k\subset \mathbb S^3$ be a
satellite
knot with $\hbox{rank }\pi_1 (E(k))=2$. Then one of the following
holds:
\begin{enumerate} \item Either $k$ has tunnel number one, or

\item
$E(k)=M_{1}\cup M_{2}$ where $M_1$ is the exterior of a $(p,q)$ torus
knot
and $M_2$ is the exterior of a hyperbolic link $L=k_{1}\cup
k_2\subset\mathbb S^3$ with two unknotted components and tunnel number
$\geq 2$.
Furthermore $\pi_1(M_2)$ is generated by two elements one of which
corresponds to a meridian $m$ and the gluing map identifies the fibre
of $M_{1}$
with $m$. \end{enumerate}
\end{Corollary}

Corollary \ref{sat} is already a consequence of the work of
S. Bleiler and A. Jones [BJ2] (cf. also [B3]).

Some of the results of this paper  have also been obtained independently by
D. Bachman,
S. Bleiler and A. Jones [BBJ], using a different combinatorial group
theoretical approach.

\section{Some combinatorial tools}

In this section we recall the results of \cite{KW} on generating
pairs of
fundamental
groups of graphs of groups and draw some conclusions. Some other
lemmas that are
needed later are also shown.

In \cite{KW} the action of 2-generated groups on simplicial trees is
investigated. Suppose that a group $G$ acts simplicially on a
simplicial tree $T$ and that $g\in G$. We define $T_g=\{x\in
T|g^zx=x\hbox{ for some }z\in\mathbb Z\hbox{ with
}g^z\neq 1\}$. It is clear that $T_g\neq\emptyset$ if and only if $g$
acts with
a fixed point. With this notation the main result of \cite{KW}
immediately
implies the following:

\begin{Theorem}\label{KW} Let $G$ be a torsion-free non-free
2-generated group acting simplicially and without inversion on the
simplicial tree $T$. Then any generating pair is Nielsen-equivalent
to a pair $\{g,h\}$ such that either

\begin{enumerate}

\item $T_g\cap T_h\neq\emptyset$ or

\item $T_g\cap hT_g\neq\emptyset$.

\end{enumerate}

\end{Theorem}

\noindent {\em Remark} It should be noted that the proof actually
shows that once
there is generating pair $\{g,h'\}$ such that $g$ acts with a fixed
point then one of the two statements hold for a pair $\{g,h\}$ where
$h=g^{z_1}h'g^{z_2}$ for some $z_1,z_2\in\mathbb Z$.

\

Following \cite{Se} we say that a splitting of a group $G$ as a
fundamental group of a
graph of groups is {\em 2-acylindrical} if no non-trivial element 
$g\in G$ fixes a
segment of length greater than $2$ in the corresponding Bass-Serre
tree. We will sometimes use the following simple fact:

\begin{Lemma} Suppose that $G$ is the fundamental group of a $2$-acylindrical graph of groups and that $T$ is the corresponding Bass-Serre tree. 

Then every vertex $v$ that is
fixed under the action of a non-trivial power of an element $g$ is in
distance at most $1$ of a vertex fixed
under the action of $g$. 
\end{Lemma}

\noindent{\em Proof} Choose $n$ such that $g^n\neq 1$ and that $g^nv=v$. Let $w$ be the vertex of $T$ that is closest to $v$ such that $gw=w$. Clearly $g^n$ fixes the segment $[v,w]$. The minimality of the distance between $v$ and $w$ guarantees that $g$ fixes no point of $[v,w]$ except $w$. It follows that the $[v,w]\cap g[v,w]=[v,w]\cap [gv,gw]=[v,w]\cap [gv,w]=\{w\}$. In particular $(v,gv)=2d(v,w)$. Now $g^n g[v,w]=g g^n[v,w]=g[v,w]$, i.e. $g$ fixes $[v,w]$ and $g[v,w]$ and therefore $[v,gv]$. As the action of $G$ is $2$-acylindrical this implies that $d(v,gv)\le 2$ and therefore $d(v,w)\le 1$ which proves the assertion.\hfill$\Box$

\medskip The following Lemma is a generalization of
Corollary~3.2 of \cite{KW}
which is in turn a generalization of the main result of \cite{BJ1}.

\begin{Lemma}\label{amalgam} Suppose that $G=A*_C B$ with $C\neq 1$
is torsion-free and that this splitting is 2-acylindrical.
If $G$ is 2-generated then there exists a generating pair $\{g,h\}$
such that (possibly after exchanging $A$ and $B$) $g\in A$ with
$g^n\in C-1$ and that one of the following holds:

\begin{enumerate}

\item $h\in B-C$.

\item $h=ab$ with $a\in A-C$, $b\in B-C$ and $a^{-1}g^ma\in C$ for
some $m\in \mathbb N$.

\item $h=bab^{-1}$ with $a^m\in C-1$ for some $m\in \mathbb N$.

\end{enumerate} \end{Lemma}

\noindent {\em Proof } Applying Theorem~\ref{KW} to an arbitrary
generating pair
guarantees the existence of a generating pair $\{g,h\}$ such that $g$
acts with a
fixed point and that either $T_g\cap T_h\neq\emptyset$ or $T_g\cap
hT_g\neq\emptyset$.

Since the splitting is 2-acylindrical we know that there
exists a vertex
$v\in T_g$ such that $gv=v$ and that for every vertex $w\in T_g$ we
have $d(v,w)\le 1$. In the case that $h$ also acts with a fixed
point, $T_h$ clearly has the same structure. After conjugation we
have that $v$ is the vertex fixed under the action of $A$ or $B$,
w.l.o.g. we assume that $Av=v$, i.e. $g\in A$.

\noindent {\bf 1. case:} $T_g\cap T_h\neq\emptyset$. Choose $z\in
T_g\cap T_h$ and $x\in T_h$ with $d(z,x)\le 1$ such that $hx=x$. Note
that $d(v,x)\le 2$. If $v=x$ then $\langle g,h\rangle\subset A$ which
contradicts our assumption that $\{g,h\}$ is a generating set of $G$.
If $d(v,x)=1$ then we can assume after conjugation with an element of
$A$ that $[x,v]$ is the edge fixed under the action of $C$ and $v$ is
the vertex fixed by $B$, in particular $h\in B$. It now follows that
either a power of $g$ or a power of $h$ must fix $[x,v]$
since otherwise $T_g\cap T_h=\emptyset$, i.e that either $g^n\in C-1$
or $h^n\in
C-1$ for some power of $g$ or $h$. If $g^n\in C-1$ then we are in
situation 1 of
the lemma, otherwise we are in situation 1 after exchanging $A$ and
$B$ and $g$ and $h$. If $d(x,v)=2$ then $[x,v]=[x,z]\cup [z,v]$ and
after conjugation with an
element of $A$ we can assume that $[v,z]$ is fixed under the action
of $C$. This
implies that $g^n\in C-1$ for some $n\in \mathbb N$. It is clear that
$x=bv$ for
some $b\in B$, it follows that $h\in bAb^{-1}$, i.e. that
$h=bab^{-1}$ for some $a\in A$. Since a non-trivial power of $h$
fixes $z$ it follows that $h^m=ba^mb^{-1}\in B-1$ for some $m\in
\mathbb N$, i.e. $a^m\in C-1$ for some $m\in \mathbb N$. This gives
situation 3.

\noindent {\bf 2. case:} $T_g\cap hT_g\neq\emptyset$. Chose $x,y\in
T_g$ such that $hx=y$. Since in the Bass-Serre tree of an amalgamated
product $G$-equivalent vertices are always is even distance we can
assume that either $d(x,y)=0$ or $d(x,y)=2$. If $d(x,y)=0$ we are in
the first case, i.e. we can assume that $d(x,y)=2$. It is clear that
$d(v,x)=d(v,y)=1$ and after conjugation
with an element of $A$ we can assume that $x$ is the vertex fixed
under the action of $B$, in particular $[v,x]$ is fixed under the
action of $C$ which implies that $g^n\in C-1$ for some $n\in\mathbb
N$. Now it is clear that $ax=y$ for some $a\in A$. This implies that
$h=ba$ for some $b\in \Stab y=B$ and $b\in B$. Since some power $g^k$
of $g$ fixes $[v,y]$ it follows that $g^{nk}$ fixes $[v,x]$ and
$[v,y]=a[v,x]$. This implies in particular that $a^{-1}g^{nk}a\in C$.
This gives situation~2. \hfill$\Box$

\medskip The next lemma gives a bound on the number of vertex groups
if the group is 2-generated and the graph underlying the splitting is
homeomorphic to a circle.

\begin{Lemma}\label{HNN1} Let $G$ be a torsion-free non-free
2-generated group. Then $G$ does not admit a 2-acylindrical
splitting whose underlying graph
is homeomorphic to a circle and has more than $2$ vertices.
\end{Lemma}

\noindent {\em Proof} Suppose that the underlying graph has at least
three vertices. We apply Proposition~\ref{KW}, i.e. we can assume
that there exists a generating pair
$\{g,h\}$ such that either $T_g\cap T_h\neq\emptyset$ or that
$T_g\cap hT_g\neq\emptyset$.

The case $T_g\cap T_h\neq\emptyset$ cannot occur since all elements
that act with a fixed point lie in
the kernel of the quotient map that quotients out all vertex groups.
In our case
however this quotient is an infinite cyclic group since there lies
one edge outside a maximal tree of the graph underlying the
splitting, a contradiction.

It remains to rule out the case $T_g\cap hT_g\neq\emptyset$, i.e.
that there exists $x,y\in T_g$ such that $hx=y$. The diameter of
$T_g$ is at most $2$ since the splitting is 2-acylindrical,
in particular $d(x,y)\le 2$. This however implies that $\langle
g,h\rangle$ lies in a vertex group of the graph of group obtained
from the original graph of groups after collapsing the at most two
edges corresponding to the edges of $[x,y]$. This however implies
that $\langle g,h\rangle\neq G$ since the resulting graph of groups
still contains at least one edge.\hfill$\Box$

\medskip The next two lemmas give some information on the generators
if there are
only two vertices and if there is only one vertex, respectively.

\begin{Lemma}\label{HNN2} Let $A$ and $B$ be two torsion-free groups
and $C_1,C_2\subset A$ and $C_3,C_4\subset B$ subgroups such that
there exist isomorphisms $\phi_1:C_1\to C_3$ and $\phi_2:C_2\to C_4$.
Let $G=(A*_{C_1=C_3}B)*_{C_2=C_4}$, i.e. $G=\langle
A,B,t|\phi_1(c_1)=c_1,\phi_2(c_2)=tc_2t^{-1}\rangle$ and suppose that
this splitting is 2-acylindrical. Then there exists a
generating pair $\{g,h\}$ such that (possibly after exchanging $A$
and $B$)

\begin{enumerate}

\item $g\in B$, $g^n\in C_3-1=C_1-1$ and $g^m\in
bC_4b^{-1}-1=btC_2t^{-1}b^{-1}$
for some $b\in B$ and $n,m\in \mathbb N$ and

\item $h=bta$ for some $a\in A$.

\end{enumerate}
\end{Lemma}

\noindent {\em Proof } We study the action on the associated
Bass-Serre tree. As
in the proof of Lemma~\ref{HNN1} we apply Theorem~\ref{KW} and
exclude the case that $g$ and $h$ act with a fixed point, i.e. we can
assume that $T_g\cap hT_g\neq\emptyset$ and that $h$ acts without
fixed point. Choose a vertex $v\in T_g$ such that $gv=v$; this
implies that $d(v,x)\le 1$ for all $x\in T_g$ since the action is 2-acylindrical. After conjugation and possibly interchanging
$A$ and $B$ we can assume that this vertex is fixed by $B$, i.e. that
$g\in B$. Choose further $x,y\in T_g$ such that $hx=y$. Note that
$d(x,y)$ is even since the segment $[x,y]$ must map onto a closed
path on the quotient graph. Since we assume that $h$ acts without
fixed point this implies that $d(x,y)=2$, i.e. that $[x,y]=[v,x]\cup
[v,y]$ where $[v,x]$ and $[v,y]$ are edges. We can assume that there
exists no $b\in B$ such that $b[v,x]=[v,y]$ since $h$ would then be a
product of elliptic elements and therefore together with $g$ lie in
the kernel of the map that quotients out the vertex
groups. It follows that we can assume that $[v,x]$ is $B$-equivalent
to the edge
associated to $C_3$ and $[v,y]$ is $B$-equivalent to the edge
associated to $C_4$. After conjugation in $B$ we can assume that
$[v,x]$ is the edge associated
to $C_3$, i.e. that $x$ is fixed under $A$. It follows that $[v,x]$
is fixed by $C_3$ and $[v,y]$ is fixed under $bC_4b^{-1}$ where $b\in
B$ is such that $btx=y$. This guarantees the first conclusion of the
lemma. Since $y=btx$ it follows that $h=sbt$ for some $s\in\Stab
y=bt(\Stab x)t^{-1}b^{-1}=btAt^{-1}b^{-1}$. This implies that
$sbt=(btat^{-1}b^{-1})bt=bta$
for some $a\in A$ which proves the Lemma. \hfill$\Box$

\begin{Lemma}\label{HNN3} Let $A$ be a group, $C_1,C_2$ be two
isomorphic subgroups with an isomorphism $\phi:C_1\to C_2$. Let
$G=A*_{C_1=C_2}=\langle A,t|tc_1t^{-1}=\phi(c_1)\rangle$ and suppose
this splitting is 2-acylindrical. Then there exists a
generating pair $\{g,h\}$ such that (possibly
after exchanging $C_1$ and $C_2$) $g\in A$ and $g^n\in C_1$ and
$h=at$ for some $a\in A$. \end{Lemma}

\noindent {\em Proof } As in the proof of Lemma~\ref{HNN2} we can
assume that there is a generating pair $\{g,h\}$ such that $T_g\cap
hT_g\neq\emptyset$ and that $h$ acts without fixed point. Chose
$x,y\in T_g$ such that $hx=y$. It is clear that $d(x,y)$ is odd since
otherwise the exponent sum of all occurrences of
$t$ in $h$ would be even and therefore not lie in $\{-1,1\}$, which
implies that
the images of $g$ and $h$ do not generate the cyclic quotient of $G$
by $N_G(A)$ which
is generated by the image of $t$. Since the diameter of $T_g$ is at
most $2$ this
implies that $d(x,y)=1$. Possibly after exchanging $x$ and $y$ we can
assume that
$gx=x$ since the action is 2-acylindrical. After conjugation
we can assume that $x$ is fixed under the action of $A$ and that
$[x,y]$ is fixed under
the action of either $C_1$ or $C_2$. Possibly after exchanging $C_1$
and $C_2$ we
can assume that $[x,y]$ is fixed by $C_1$, i.e. that $y=t^{-1}x$. It
follows that
$h=st^{-1}$ where $s\in \Stab y=t^{-1}At$, i.e.
$h=(t^{-1}at)t^{-1}=t^{-1}a$ for
some $a\in A$. After replacing $h$ with its inverse this proves the
lemma.\hfill$\Box$

\medskip

In some instances it is important to see that a set $S$ does not
generate a given
group $G$. If $G$ is given as a fundamental group of a graph of
groups this can be seen if the induced splitting of the subgroup
generated by $S$ is distinct from the
original splitting of $G$. In \cite{FRW} situations where investigated
when the induced splitting can be read of a particular generating
set. The following Proposition describes the situations needed in the
course of this paper. It is a direct consequence of the discussion in
\cite{FRW}.

\begin{Proposition}\label{controlled} Let $G$ be a group that acts
simplicially without inversion on a simplicial tree $T$. Let $ET$ be
the set of edges of $T$.

\begin{enumerate}

\item Let $e=[v,w]\in ET$. Suppose that $G_v,G_w\subset G$ such that
$G_vv=v$, $G_ww=w$ and $G_v\neq G_v\cap\Stab e=G_w\cap \Stab e\neq
G_w$. Then the induced splitting of $U=\langle G_v,G_w\rangle$ is
$G_v*_{G_v\cap\Stab e}G_w$.

\item Let $e_1=[v,w],e_2=[v,z]\in ET$. Suppose that
$G_v,G_w,G_z\subset G$ such that $G_vv=v$, $G_ww=w$ and $G_zz=z$.
Suppose further that $G_v\cap\Stab e_1=G_w\cap \Stab e_1\neq G_w$,
that $G_v\cap\Stab e_2=G_z\cap \Stab e_2\neq G_z$ and that $e_1$ and
$e_2$ are not $G_v$ equivalent. Then the induced splitting of
$U=\langle G_v,G_w,G_z\rangle$ is $G_w*_{G_v\cap\Stab
e_1}G_v*_{G_v\cap\Stab e_2}G_z$.

\item Let $e=[v,w]\in ET$. Suppose that $G_v,G_w\subset G$ and $h\in
G$ such that $G_vv=v$, $G_ww=w$, $G_v\cap \Stab e=G_w\cap\Stab e$,
$hv=w$ and $h^{-1}G_wh=G_v$. Suppose further that $U=\langle
G_v,h\rangle$. Then
the induced splitting of $U$ is the HNN-extension of $G_v$ along
$G_v\cap \Stab e$.

\item Let $e_1=[v,w],e_2=[v,z]\in ET$. Suppose that
$G_v,G_w,G_z\subset G$ and $h\in G$ such that $G_vv=v$, $G_ww=w$,
$G_zz=z$, $G_v\cap \Stab e_1=G_w\cap\Stab e_1$, $G_v\cap\Stab
e_2=G_z\cap\Stab e_2$, $w=hz$, $G_z=hG_wh^{-1}$, that $e_1$ and $e_2$
are not $G_v$-equivalent and that $e_1$ and $he_2$ are not
$G_w$-equivalent. The the induced splitting of $U=\langle
G_v,G_w,h\rangle$ is $(G_v*_{C_1}G_w)*_{C_2}$ where $C_1=G_v\cap \Stab
e_1=G_w\cap\Stab e_1$ and $C_2=G_v\cap \Stab e_2=G_z\cap\Stab
e_2$.~\hfill$\Box$

\end{enumerate}

\end{Proposition}

\section{Some topological lemmas}\label{toplem}

If $M^3$ is a Seifert manifold with base orbifold $\B$ we denote
$\pi_1(M^3)$ by
$G$, the element of $G$ corresponding to the fibre by $f$ and
$\pi_1(\B)$ by $F$.
We then have the following exact sequence $$1\to \langle
f|-\rangle\to G\overset{\pi}{\to} F\to 1.$$ The subgroups of $F$
corresponding to boundary curves of the base orbifold are the images
of the subgroups of $G$ corresponding
to the boundary components of $M^3$.

The only Seifert manifolds with boundary whose fibration is not
unique up to isotopy are $\mathbb R\times S^1\times S^1$ and the
orientable circle bundle over the M\"obius band, denoted as $Q$ by
Waldhausen \cite{Wa1},
which can also be fibered over $D(2,2)$. If $Q$ is fibered over the
M\"obius band
we denote the fibre by $f_Q$ and if $Q$ is fibered over $D(2,2)$, we
denote the
fibre by $f'_Q$.

Considering $Q^3$ as the orientable circle bundle over the M\"obius
band yields the presentation $\pi_1(Q^{3})=\langle f_{Q},s
|sf_{Q}s^{-1}=f_{Q}^{-1}\rangle$,
the standard presentation of the Klein bottle group. Considering it
as the Seifert
manifold over $D(2,2)$ we get the presentation $\pi_1(Q^3)=\langle
x,y,f_Q'|[x,f_Q'],[y,f_Q'],x^2=f_Q',y^2=f_Q'\rangle$. The isomorphism
is given by
$f_{Q}=xy{f_{Q}'}^{-1}$ and $s=x$.

Let $C$ be the normal subgroup of $\pi_1(Q^{3})$ corresponding to the
boundary which is generated by $\{s^2,f_{Q}\}$ and $\{xy,f_{Q}'\}$
respectively. We see that for any $g\in \pi_1(Q^{3})-C$ we have that
$g^2\in \langle f_Q'\rangle$ and
$g(xy)^n{f_Q'}^mg^{-1}=x^{-1}(xy)^n{f_Q'}^mx=(xy)^{-n}{f_Q'}^{m+2n}$.

\begin{Lemma}\label{top1} Let $M^3$ be a $\partial$-incompressible
Seifert-manifold, $T$ a boundary component and $P\subset G$ be a
subgroup corresponding to $T$. Suppose further that $g\in G$
generates a maximal cyclic subgroup such that $\langle g\rangle\cap
P=\langle g^n\rangle$ for some $n\in \mathbb N$ and that $P_1\subset
P$ with $g^n\in P_1$. Then $\langle P_1,g\rangle
\cap P=P_1$. \end{Lemma}

\noindent {\em Proof }If $g\in P$ the assertion is trivial, i.e. we
can assume that $g\in G-P$. If $M^3$ is $Q^3$ we write it as the
Seifert space over $D(2,2)$. Let now $w$ be the element of $F$ that
corresponds to the boundary curve of $\B$ that comes from $T$. Since
$\B$ is not the M\"obius band we know that
all roots of $w^n$ are of type $w^k$. This implies that for any
element $g\in G$
such that $\pi(g^n)\in \langle w\rangle -1$ for some $n\ge 1$ we have
$g\in P$. It follows that $\pi(g^n)=1$, i.e. that $g^n\in \langle
f\rangle$. Since we assume
that $g$ generates a maximal cyclic subgroup we have $g^n=f^{\pm 1}$;
possibly after replacing $g$ with its inverse we can assume that
$g^n=f$.

The fundamental group $F$ of the base space is a free products of
cyclic groups,
i.e. $F=\langle s_1,\ldots ,s_r|s_1^{n_1},\ldots ,s_r^{n_r}\rangle$
with $n_i\in
\mathbb N\cup \{\infty\}$. It is clear that $P_1=\langle
f,m^k\rangle$ for some $k\in\mathbb N$ where $m\in P$ is chosen such
that $P=\langle f,m\rangle$. Choose
now $l$ such that $\langle P_1,g\rangle \cap P=\langle f,m^l\rangle$;
we have to
show that $k=l$. The projection $\pi$ is an injection when restricted
to $\langle m\rangle $, has kernel $\langle f\rangle$, maps $m$ onto
$w$ (or $w^{-1}$) and $g$ onto a conjugate of $s_i$ for some $i\in
\{1,\ldots,r\}$. This
means that $\pi(P_1)=\langle\pi(m^k)\rangle=\langle w^{\pm k}\rangle$
and $\pi(\langle P_1,g\rangle)\cap \langle \pi(m)\rangle=\langle
\pi(m^l)\rangle=\langle w^{\pm l}\rangle$. The assertion now follows
from the fact that $\pi(\langle P_1,g\rangle)\cap \langle
\pi(m)\rangle=\langle w^k,s\rangle\cap \langle w\rangle$ for some
$s\in F$ that is conjugate to some $s_{i}$ and that in $F$ the
statement $\langle w^k,s\rangle\cap \langle w\rangle=\langle
w^k\rangle$ holds for any element $s$ that is conjugate to one of the
$s_i$. This can be seen by looking at the orbifold group obtained by
adding the relation $w^k$ to $F$: The resulting orbifold is not bad
since $\B$ was not a disk with less than two cone points, and the
image of $s$ in the quotient is not conjugate to an elliptic element
corresponding to the new cone point.~\hfill$\Box$

\begin{Lemma}\label{top2} Let $M^3$ be a $\partial$-incompressible
Seifert-manifold with boundary component $T$ and $P\subset G$ be a
corresponding
subgroup. Suppose further that $g\in G$ such that $g^n\in P$ for some
$n\in \mathbb N$ and that $G=\langle P,g\rangle$. Then $M^3$ is a
Seifert manifold with
basis $D(p,q)$ or $A(p)$ and $g$ is a root of the fibre. \end{Lemma}

\noindent {\em Proof } The fact that $g$ is a root of the fibre
follows as in the
proof of Lemma~\ref{top1}, possibly after rewriting $Q^3$ as the
Seifert space with base space $D(2,2)$. Now the base space must have
at least one boundary component and be generated by the element
corresponding to the boundary and a torsion element. The only base
spaces with this property are $D(p,q)$ and $A(p)$.
\hfill$\Box$

\begin{Lemma}\label{top3} Let $M^3$ be an orientable Seifert-manifold
which is not $T^2\times I$. Suppose that $T_1$ and $T_2$ are two
boundary components with
corresponding subgroups $P_1$ and $P_2$. Suppose that $g\in P_1$ and
$h\in P_2$,
that neither $g$ nor $h$ correspond to the fibre and that $\langle
g,h\rangle$ is
not free. Then the following hold:

\begin{enumerate}

\item $\B$ is of type $A(p)$, i.e. $G=\langle s,x,f|[s,f],\allowbreak
[x,f],x^p=f^b\rangle$ with $1\le b<|p|/2$ and $(p,b)=1$ and after
conjugation we
have $g=sf^m$ and $h=xsf^n$ for some $m,n\in \mathbb Z$. In
particular $\langle g,h\rangle$ maps surjectively onto $F$.

\item If $\langle g,h\rangle =G$ then we have additionally that $b=1$
and after conjugation we have $g=sf^m$ and $h=xsf^m$ for some $m\in
\mathbb Z$. In particular $M^3$ is the exterior of a 2-bridge link
and $g$ and $h$ correspond to
the meridians.

\end{enumerate} \end{Lemma}

\noindent {\em Proof } We first show that $H:=\pi(\langle g,h\rangle
)=\langle \pi(g),\pi(h)\rangle\subset F$ is free in $\pi(g)$ and
$\pi(h)$ unless $\B$ is of type $A(p)$ and $\pi(g)$ and $\pi(h)$
correspond to the
boundary curves of $\B$. This clearly implies the first assertion of
the lemma.

Note that $\pi(g)$ and $\pi(h)$ are powers of elements that
correspond to the boundary curves of $\B$. Since $H$ is a subgroup
of a free product of cyclic group we know by Kurosh's subgroup
theorem that $H$ itself is a free product of cyclic groups. It
follows that $H$ is free if and only if $H$ is
torsion-free.

Suppose now that either $\B$ is not of type $A(p)$ or that $\B$ is of
type $A(p)$
and $\pi(g)$ (the case of $\pi(h)$ is analogous) is a proper power of
a boundary
curve. We look at the quotient map $\phi:F\to F/N_F(g)$. It is clear
that no torsion element lies in the kernel of $\phi$ since $\B$ has
at least two boundary
components and the resulting orbifold is therefore good. In
particular $\phi (H)$ has torsion if $H$ has torsion. It is clear
that $\phi (H)=\langle \phi(\pi(h))\rangle$. The group $\langle
\phi(\pi(h))\rangle$ however is infinite cyclic since it is a
subgroup of the cyclic subgroup of an orbifold group corresponding to
a boundary curve which is infinite because the orbifold is different
from $D$ (since $\B$ was assumed to not
be of type $A$) and $D(p)$ (note, that if $\B$ was of type $A(p)$ the
new orbifold
is of type $D(p,q)$ since we assumed that $g$ was a proper power of
the boundary
curve). It follows that $H$ is torsion free and therefore free.

It follows that either $\langle g,h\rangle$ is free or that
$G=\langle s,x,f|[s,f],\allowbreak [x,f],x^p=f^b\rangle$ for some
$p>b>0$ with $b<|p/2|$ and
$(b,p)=1$ and that $g$ is conjugate to an element of type $sf^k$ and
$h$ is conjugate to an element of type $(xs)f^l$ for some
$k,l\in\mathbb Z$.

\noindent {\bf Claim: }Either $\langle \pi(g),\pi(h)\rangle$ is free
in $\pi(g)$
and $\pi(h)$ or after conjugation $\pi(g)=s$ and $\pi(h)=xs$.

\noindent After a suitable conjugation $\pi(g)$ is in the desired
form, i.e. $\pi(g)=s$ and $\pi(h)=w(xs)w^{-1}$ for some $w\in
F=\langle s,x|x^p\rangle=\langle s|-\rangle *\langle x|x^p\rangle$.
It is easy to see that
$\pi(h)$ has normal form of one of the types $y_1\cdots
(y_{l}x)sy_{l}^{-1}\cdots
y_{1}^{-1}$, $y_1\cdots (y_{l}s)xy_{l}^{-1}\cdots y_{1}^{-1}$,
$y_1\cdots y_{l}x(sy_{l}^{-1})\cdots y_{1}^{-1}$ or $y_1\cdots
y_{l}s(xy_{l}^{-1})\cdots y_{1}^{-1}$. If $y_{1}\in\langle s\rangle$
we can conjugate $\pi(g)$ and $\pi(h)$
by $y_{1}$. This conjugation does not change $\pi(g)$ but reduces the
length of $\pi(h)$; i.e. we can assume that $y_{1}\in\langle
x\rangle$. If $l\ge 2$ we see
that no cancellation occurs in products of powers of $\pi(g)$ and
$\pi(h)$, since
both $\pi(g)$ and $\pi(h)$ are of infinite order this implies that
$\langle \pi(g),\pi(h)\rangle$ is free in $\pi(g)$ and $\pi(h)$. If
$l=1$ cancellation occurs if and only if $y_{1}=y_{l}=1$. It follows
that either $\pi(h)=xs$ or $\pi(h)=sx$.
After conjugation of the pair $\{g,h\}$ we have $\pi(h)=xs$ which
proves the claim.

If $\langle g,h\rangle$ is not free in $g$ and $h$ it follows that
$\langle \pi(g),\pi(h)\rangle$ is not free in $\pi(g)$ and $\pi(h)$,
i.e. after conjugation $g=sf^k$ and $h=(xs)f^l$. It is clear that
$\langle g,h\rangle=\langle sf^k,(xs)f^l\rangle=\langle
sf^k,xf^{l-k}\rangle$ maps surjectively to $F$, it remains to verify
the second assertion, i.e. we
have to determine the situations where in addition $\langle
g,h\rangle\cap \langle f\rangle=\langle f\rangle$. Now a freely
reduces products in $x$ and $s$ is trivial in $F$ if and only if it
is a product of conjugates of $p_{th}$ powers of $x$. It follows that
a freely reduced product in $sf^k$ and $xf^{l-k}$ lies in the kernel
of $\pi$ if and only if it is a product of conjugates of $p_{th}$
powers of $xf^{l-k}$. Since
$w(xf^{l-k})^pw^{-1}= wx^pf^{p(l-k)}w^{-1}=
wf^bf^{p(l-k)}w^{-1}=wf^{b+a(l-k)}w^{-1}=
f^{b+a(l-k)}$ for any $w\in F$ it follows that $\langle
g,h\rangle\cap\langle
f\rangle=\langle f^{b+a(l-k)}\rangle$. This implies that
$b+a(l-k)=\pm 1$ and therefore $b=1$ and $l=k$.~\hfill$\Box$

\begin{Lemma}\label{top5} Let $M^3$ be an orientable
$\partial$-incompressible Seifert-manifold which is not $T^2\times
I$. Suppose that $T$ is a boundary component. Let $P\subset
\pi_1(M^3)$ be the corresponding subgroup and $g\in P$ where $g$ is
primitive in $P$ and does not correspond to the fibre. Then there
exists an element $h\in\pi_1(M^3)$ such that $\langle
g,h\rangle=\pi_1(M^3)$ if and only if
one of the following holds (after conjugation):

\begin{enumerate}

\item $\B$ is of type $A(p)$, i.e. $G=\langle s,x,f|[s,f],\allowbreak
[x,s],x^p=f^b\rangle$ and $g=sf^k$ or $g=sxf^k$ for some $k$.

\item $\B$ is of type $D(p,q)$, i.e. $G=\langle
x,y,f|[x,f],[y,f],x^{p}=f^{b_1},y^{q}=f^{b_2}\rangle$ and $g=xyf^n$
for some $n\in\mathbb N$.

\item $\B$ is of type $M\ddot o(p)$, i.e. $G=\langle
x,s,f|[x,f],x^p=f^b,sfs^{-1}=f^{-1}\rangle$ and either $b=1$ and
$g=s^2x$ or $a=1$, $b=1$ and $g=s^2xf^{-1}$.

\item $\B$ is of type $M\ddot o$, i.e. $G=\langle
s,f|sfs^{-1}=f^{-1}\rangle$ and
$g=s^{2l}f^{\pm 1}$ for some $l\in\mathbb Z$.

\end{enumerate}

\noindent Except in the first case this implies that $M^3$ is the
exterior of a 1-bridge knot in a lens space and $g$ corresponds to
the meridian. \end{Lemma}

\noindent {\em Proof } Let $M(g)$ be the manifold obtained by a Den
filling of $M$ along the (simple) curve on $T$ corresponding to $g$.
Since $g$ is not a power of the fibre we can extend the Seifert
fibration of $M$ to $M(g)$. The classification of small Seifert
manifolds shows that $M(g)$ has cyclic fundamental group, i.e. is a
lens space, ifandonlyif $M$ and $g$ are as in one of the cases
(1)-(4) of
Lemma~\ref{top5}. Except in the first cases this implies that $M^3$
is the exterior of a 1-bridge knot in a lens space and $g$
corresponds to the meridian, see Lemma~1 of \cite{W2}. This
guarantees in particular the existence of the appropriate $h$. In
the first case $h$ can be chosen as a generator of the cyclic
subgroup $\langle x,f\rangle$.\hfill$\Box$

\begin{Lemma}\label{top10} Let $M^3$ be an orientable
$\partial$-incompressible Seifert-manifold that is not $T^2\times I$
with a boundary component $T$ and $P\subset G=\pi_1(M^3)$ be a
corresponding subgroup.

Suppose further that $g\in P$ and $h\in
wPw^{-1}-P$ for some $w\in G$ and that $g$ does not correspond to a
power of the
fibre. Denote the intersection number of $g$ and $h$ with the fibre
on $T$ by $n_g$ and $n_h$.

Then $\langle g,h\rangle$ is free in $g$ and $h$ or $M^3$ is a
Seifert space with base space of type $D(2,q)$, $\min
(|n_g|,|n_h|)=1$ and $\max(|n_g|,|n_h|)\le 3$. If $\max
(|n_g|,|n_h|)>1$ or $q$ is even we further have that $h$ is not
conjugate
to an element of $P$ in $\langle g,h\rangle $. \end{Lemma}

\noindent {\em Proof } Possibly after conjugation and exchanging $g$
and $h$ and
replacing $g$ or $h$ by their inverses we can assume that $1\le
n_g\le n_h$ since
neither $g$ nor $h$ corresponds to a power of the fibre. We can
assume that $M^3$
is not $Q^3$ since in this case $wPw^{-1}=P$ for all $w\in G$. We
first show that $\langle g,h\rangle$ is free in $g$ and $h$ unless the
base $\B$ is of type $D(p,q)$. We actually show that the projections
$\pi (g)$ and $\pi(h)$ generate a free subgroup in the base group.
Note first that $\langle \pi(g),\pi(h)\rangle$ is not cyclic since
the cyclic subgroup corresponding to the boundary curve is malnormal
and therefore the unique maximal cyclic subgroup that contains
$\pi(g)$. However $\pi(h)$ does not lie in this subgroup. Since the
base group is a free product of cyclic it follows from Kurosh's
theorem that $\langle \pi(g),\pi(h)\rangle$ is the free product of
cyclics. To conclude we need to show that $\langle
\pi(g),\pi(h)\rangle$ is torsion-free unless $\B$ is of type
$D(p,q)$. It is clear that $\langle \pi(g),\pi(h)\rangle$ lies in the
kernel of
the map that quotients out the element corresponding to the boundary
curve of $\B$
that corresponds to $T$. This kernel however only contains torsion
elements if the resulting orbifold is bad, i.e. if $\B$ was of type
$D$, $D(p)$ or $D(p,q)$. The first two cases cannot occur since we
assume $M^3$ to be $\partial$-incompressible, i.e. the claim holds.

Suppose now that $\B$ is of
type $D(p,q)$. The base group $F$ has the presentation $\langle
a,b|a^p,b^q\rangle=\langle a|a^p\rangle*\langle b|b^q\rangle$ and
$\g=(ab)^{n_g}$
and $\h=v(ab)^{n_h}v^{-1}$ for some $v\in F$. We write $v=x_1\cdots
x_k$ as a normal form with respect to the free product $\langle
a|a^p\rangle*\langle b|b^q\rangle$. It follows that $\h=x_1\cdots
x_k(ab)^{n_h}x_k^{-1}\cdots x_1^{-1}$. W.l.o.g. we can assume that
the free product length of $\h$ is minimal
with respect to conjugation with powers of $ab$ since this
corresponds to conjugation of the pair $\{\g,\h\}$. The normal form
of $\h$ is clearly of one of
the types $x_1\cdots (x_la)b(ab)^{n_h-1}x_l^{-1}\cdots x_1^{-1}$ with
$x_la\in\langle a\rangle -1$, $x_1\cdots
x_l(ab)^{n_h-1}a(bx_l^{-1})\cdots x_1^{-1}$ with $bx_l^{-1}\in\langle
b\rangle -1$, $x_1\cdots (x_lb)a(ba)^{n_h-1}x_l^{-1}\cdots x_1^{-1}$
with $x_lb\in\langle b\rangle -1$ or
$x_1\cdots x_l(ba)^{n_h-1}b(ax_l^{-1})\cdots x_1^{-1}$ with
$ax_l^{-1}\in\langle
a\rangle -1$.

If $l\ge 2$ then $x_1\neq b^{-1}$ and $x_1\neq a$ since otherwise we
could reduce
the length of $\h$ by conjugation with $ab$ or $(ab)^{-1}$. It
follows that no cancellation
occurs in products in $\g$ and $\h$ with implies that $\langle
\g,\h\rangle $ is
free in $\g$ and $\h$.

Suppose now that $l=1$. We carry out the case
$\h=(x_1a)b(ab)^{n_h-1}x_1^{-1}$, the other cases are analogous. If
no cancellation occurs we argue as before, i.e.
we only have to study the cases (i) $x_1a=a$, i.e. $x_1=1$ and (ii)
$x_1=a$. In the first case this gives a contradiction to the
assumption that $\h=(x_1a)b(ab)^{n_h-1}x_1^{-1}$ is a normal form. In
the second case we have
$\h=a^2(ba)^{n_h-1}ba^{-1}$ and after conjugation of $\{\g,\h\}$ with
$a$ we have
$\g=(ba)^{n_g}$ and $\h=(ab)^{n_h}$. Since we assume that $(p,q)\neq
(2,2)$ we easily see that in a product of length two in $\pi(g)$ and
$\pi(h)$ at most one letter cancels. This shows that $\langle
\pi(g),\pi(h)\rangle$ is free in $\pi(g)$ and $\pi(h)$ if
$\min(n_g,n_h)\ge 2$, i.e. if the length of a reduced form is at
least $4$.

Suppose that $n_g=1$, i.e. $\g=ab$ and $\h=(ba)^{n_h}$. If $p\neq 2$
and $q\neq 2$ we see as before that $\langle \g,\h\rangle$ is free in
$\g$ and $\h$ since no cancellation occurs.

Suppose that $p=2$, the case $q=2$ is analogous. It follows that
$\langle \g,\h\rangle=\langle ab,(ba)^{n_h}\rangle =\langle
ab,(ba)^{n_h}\cdot ab\rangle =\langle ab,(ba)^{n_h-1}b^2\rangle$. If
$n_h\ge 4$ or $q\ge 4$ and $|n_h|\ge 2$ then $(ba)^{n_h-1}b^2$ is of
infinite order and any power has normal form that starts with $b$ and
ends with $b^2$. Again no cancellation occurs and we see that
$\langle \g,\h\rangle$ is free.

If $q=3$ and $n_h\in\{2,3\}$ or $q\ge 4$ even and $n_h=1$ then
cancellation arguments show that $\langle \g,\h\rangle=\langle
ab,(ba)^{n_h-1}b^2\rangle =\langle ab\rangle *\langle
(ba)^{n_h-1}b^2\rangle\cong \mathbb Z*\mathbb Z_k$ for some
$k\in\mathbb Z$ (which depends on the situation we are in). The
normal form of $\h$ with respect to this free product has length $2$
and is clearly not conjugate to an element of $\langle ab\rangle$,
which have normal form of length $1$. This implies that $h$ is not
conjugate to an
element of $P$ in $\langle g,h\rangle$.\hfill$\Box$

\begin{Lemma}\label{top6} Let $M^3$ be an orientable
$\partial$-incompressible Seifert-manifold which is not $T^2\times I$
and not $Q^3$ with a boundary component $T$ and $P\subset G=
\pi_1(M^3)$ be a corresponding subgroup. Let $g\in P$ such that $g$
is not the fibre. Suppose that $h\in \pi_1(M^3)-P$ such that $\langle
g,hgh^{-1}\rangle$ is not free. Then then one of the following holds:

\begin{enumerate}

\item If $\langle g,hgh^{-1}\rangle=G$ for some $h\in G$ then $M^3$
is the exterior of the
(2,p)-torus knot and $g$ corresponds to a meridian. In particular
$M^3$ is the exterior of a 2-bridge knot.

\item If $\langle g,hgh^{-1}\rangle$ maps surjectively onto the base
group for some $h\in G$ then the base orbifold is of type $D(2,p)$
with odd $p$ and $g$ maps onto the element of the base group
corresponding to a boundary curve.

\item The base manifold is of type $D(2,2l)$ and $g$ maps onto the
element of the base group corresponding to a boundary curve. In this
case $hgh^{-1}$ is not conjugate to an element of $P$ in $\langle
g,hgh^{-1}\rangle$.

\end{enumerate}
\end{Lemma}

\noindent {\em Proof } It is clear that $g$ and $hgh^{-1}$ have the
same non-trivial intersection number $n$ with the fibre. It follows
therefore immediately from Lemma~\ref{top10} that the base is of type
$D(2,n)$. It further follows from the proof of Lemma~\ref{top10} that
we can assume that $\pi(g)=ab$ and $\pi(h)=ba$. If $n$ is odd this
implies that $\pi(g)$ and $\pi(h)$ generate $F$ which puts us into
situation 2. If $n$ is even Lemma~\ref{top10} implies that we are in
situation 3.

To show that the first statement holds we look at the manifold
$M^3(g)$ obtained from $M^3$ by a Dehn filling killing the curve
corresponding to $g$. Since $g$ has intersection number one with the
fibre, we can extend the Seifert fibration of $M^3$ to a Seifert
fibration of $M^3(g)$. Since $M^3(g)$ has trivial fundamental group
it must be $\mathbb S^3$. It follows that $M^3$ is a Seifert fibered
knot exterior in $\mathbb S^3$ that has base space $D(2,q)$. This
implies that $M^3$ is the exterior of the (2,p)-torus knot. Now $g$
must correspond to a meridian since torus knots have property
P.~\hfill$\Box$

\medskip

\begin{Proposition}\label{twobridge} Let $M^3$ be a compact
orientable 3-manifold
with a complete hyperbolic structure of finite volume on its
interior. Suppose that $U$ is a subgroup of $\pi_1(M^3)$ which is
generated by two parabolic primitive elements. Suppose furthermore
that these two parabolic elements are conjugated in $\pi_1(M^3)$ if
$\partial M^3$ is connected.

Then either $U$ is
free or $U$ is Abelian or one of the following holds:

\begin{enumerate}
\item $U=\pi_1(M^3)$ and $M^3$ is the exterior of a 2-bridge knot or
link in $\mathbb S^3$.

\item $|\pi_1(M^3):U|=2$ and the covering space $\hat M^3$ of $M^3$
corresponding to $U$ is the exterior of a 2-bridge link in $\mathbb
S^3$ (with 2 components).
\end{enumerate}
Moreover the two parabolic generators of $U$ correspond to meridian
curves of the 2-bridge knot or link.
\end{Proposition}

\noindent {\em Proof} Assume that $U$ is neither Abelian, nor free.
Since $M^3$ is irreducible and atoroidal, by [JS2, Thm.VI.4.1] $U$
must be of finite index in
$\pi_1(H^3)$. Hence $U \cong \pi_1(\hat M^3)$, where $\hat M^3$ is a
finite covering of $M^3$. In particular the interior of $\hat M^3$
admits a complete hyperbolic structure with finite volume. The proof
of Proposition~\ref{twobridge}
follows now from Lemma~\ref{hyp1}, Lemma~\ref{hyp2} and the work of
M.~Sakuma on symmetries of spherical Montesinos links \cite{Sa2}.

\begin{Lemma}\label{hyp1} $\hat M^3$ is homeomorphic to the exterior
of a 2-bridge knot or link $L \subset \mathbb{S}^3$ and the two
parabolic
generators correspond to meridians of the knot or the link
$L$.\end{Lemma}

\noindent {\em Proof} Since $\pi_1(\hat M^3)$ is generated by two
parabolic elements, the proof follows essentially from [Ad, Thm.3.3],
together with [BZm, Prop.3.2]. We show here how to use Thurston's
orbifold theorem (cf.[BP], [CHK]) to avoid
to assume Poincar\'e conjecture in the proof.

A homological argument shows that $\partial \hat M^3$ has at most two
components. These must be tori, since $\hat M^3$ admits a complete
hyperbolic metric of finite volume on its interior. According to [Ad,
Corollary 3.2], the two generators are conjugate ifandonlyif $\hat
M^3$ has
only one torus boundary component. Moreover by [Ad, Thm.2.2] each
generator corresponds (up to conjugation) to a simple loop on
$\partial \hat M^3$. By gluing a solid torus to
each boundary component so that a meridian of the solid torus goes to
the boundary curve corresponding to a parabolic generator, one
obtains a homotopy
sphere $\Sigma^3$, that may or may not be irreducible.

Hence $\hat M^3$ is the exterior of a hyperbolic knot or link $L$ in
the homotopy
sphere $\Sigma^3$. Moreover, $\pi_1(\Sigma^3 - L)$ is generated by
two meridians.
To show that $\Sigma^3$ is in fact the true sphere $\mathbb{S}^3$ and
$L$ is a 2-bridge knot or link, we follows essentially the arguments
in [BZ, Prop.3.2].

Let $V^3$ be the 2-fold covering of $\Sigma^3$ branched along $L$.
then one has the exact sequence: $\{1\} \to \pi_1(V) \to
\pi_1(\Sigma^3 -L)/N \to \mathbb{Z}_2
\to \{1\}$, where $N$ is the subgroup of $\pi_1(\Sigma^3 -L)$
normally generated
by the squares of all meridians of $L$. Then the group
$\pi_1(\Sigma^3 -L)/N$ is
a dihedral group or a cyclic group of order 2. In the last case, by
the proof of
the Smith conjecture, $L$ would be a trivial knot contradicting that
it is a hyperbolic knot or link. Therefore $\pi_1(V)$ must be cyclic.

Since $L$ is a hyperbolic knot or link, $\Sigma^3 - L$ is irreducible
and does not contain any essential properly embedded annulus. Hence
by the equivariant sphere theorem [DD], $V$ is irreducible and
$\pi_1(V)$ is finite cyclic. By Thurston's orbifold theorem
(cf.[BP],[CHK]), $V$ is geometric, hence it is a lens space. Moreover the
covering involution is conjugated to an isometry of the spherical
structure on $V$. Hence the quotient $\Sigma^3$ is the true sphere
$\mathbb S^3$ and the branching set $L$ is a 2-bridge knot or link.
\hfill$\Box$

We use now that the two parabolic generators of $U = \pi(\hat M^3)$
are primitive in $\pi_1(M^3)$, and conjugated in $\pi_1(M^3)$ if
$\partial M^3$ is connected, to show:

\begin{Lemma}\label{hyp2} The finite covering $p: \hat M^3 \to M^3$
is regular.\end{Lemma}

\noindent {\em Proof} Denote the two parabolic generators by $g$ and
$h$. By Lemma \ref {hyp1} $\hat M^3$ is the exterior of a 2-bridge
knot
or link $L \subset \Bbb S^3$ such that the two parabolic generators
$g$ and $h$ correspond to simple closed meridian
curves of  $L$ on $\partial \hat M^3$.

If $\hat M^3$ is the exterior of a 2-bridge knot we choose a meridian
curve
$\hat \mu \subset \partial \hat M^3$ that corresponds to the free
homotopy class of $g$ and $h$.
If $\hat M^3$ is the exterior of a 2-bridge link then we choose two
meridian curves
$\hat\mu_1$ and $\hat\mu_2$ on the different components of
$\partial\hat M^3$ that represent the free homotopy classes of $g$
and $h$.

Since $g$ and $h$ are primitive in $\pi_1(M^3)$, analogously we can
choose either a
closed simple curve
$\mu$ or closed simple curves $\mu_1$ and $\mu_2$ on
$\partial M^3$, depending on whether $\partial M^3$ is connected,
which represent the
free homotopy classes of $g$ and $h$ in $\pi_1(M^3)$.

The choice of the curves clearly guarantees that either $p(\hat\mu)$
if parallel to $\mu$ on $\partial M^3$, or that $p(\hat\mu_1)$ and
$p(\hat\mu_2)$ are parallel to $\mu$ or that $p(\hat\mu_1)$ is
parallel to $\mu_1$ and that $p(\hat\mu_2)$ is parallel to $\mu_2$.

Moreover the fact that $g$ and $h$ are primitive in $\pi_1(M^3)$
implies that in any case
each component of $p^{-1}(\hat \mu))$ or $p^{-1}(\hat \mu_i)), i \in
\{1,2\}$ is mapped
homeomorphically under $p$ on to $p(\hat\mu)$ or $p(\hat\mu_i), i \in
\{1,2\}$.

Hence, the covering map $p: \hat M^3 \to M^3$ extends to a true
covering map $\bar p: \mathbb S^3 \to M^3(\mu)$ (or $\bar p: \mathbb
S^3 \to M^3(\mu_1,\mu_2)$), where $M^3(\mu)$ ($M^3(\mu_1,\mu_2)$) is
the closed orientable
3-manifold obtained by gluing a solid torus (two solid tori) to
$\partial M^3$, so that a
meridian of the solid torus goes to the boundary curve $\mu$ (the
meridians of the solid tori go to $\mu_1$ and $\mu_2$). Such a
covering
$\bar p$ is regular, hence the covering $p$ is regular.\hfill$\Box$
\medskip

We need now the following lemma:

\begin{Lemma}\label{hyp3} Any orientation preserving finite order
symmetry of $\hat H^3$, without fixed point on $\partial \hat M^3$
extends to a finite order
symmetry of $\mathbb S^3$ preserving a 2-bridge knot or
link.\end{Lemma}

\noindent {\em Proof} Since any two bridge link admits an order two
symmetry that
exchanges its components, it suffices to prove the lemma when the
symmetry preserves each component of $\partial \hat M^3$. Let
$L\subset \mathbb S^{3}$ be
a 2-bridge knot or link whose exterior is $\hat M^3$. Let $f:\hat H^3
\to \hat M^3$ be a free finite order orientation preserving symmetry,
that preserves each
component of $\partial \hat M^3$.

A free orientation preserving diffeomorphism of order $n$ of a torus
$T^2 = S^1 \times S^1$ is conjugated by a diffeomorphism isotopic to
the identity to a map of the form: $$ g(\theta,\phi) = (\theta + 2\pi
r/p, \phi + 2\pi s/q),$$ where $(r,p) = (s,q) =1$ and $lcm(p,q) = n$.
Hence it is isotopic to the identity on $T^2$.

It follows that the restriction of $f$ to each component of $\partial
\hat M^3$
is isotopic to the identity. Hence it preserves the isotopy class of
the meridian
curves of $L$ on each component of $\partial \hat M^3$. Thus, $f$
extends to a symmetry
of $\mathbb S^3$ preserving $L$. \hfill$\Box$ \medskip

\begin{Remark} Lemma~\ref{hyp2} and Lemma~\ref{hyp3} show
that, when $\partial M^3$ is connected, but $\partial \hat M^3$ has
two components, $p: \hat
M^3 \to M^3$ is a regular covering ifandonlyif the two primitive
parabolic
generators of
$\pi_(\hat M^3)$ are conjugate in $\pi_1(M^3)$.
\end{Remark}

\medskip

The following lemma finishes the proof of
Proposition~\ref{twobridge}.

\begin{Lemma}\label{hyp4} The finite covering $p: \hat M^3 \to M^3$
is either 1-sheeted or 2-sheeted in the case that $\hat M^3$ is the
exterior of a 2-bridge
link $L \mathbb S^3$. \end{Lemma}

\noindent {\em Proof}
Since the finite covering $p: \hat M^3 \to M^3$ is regular, the
finite group of covering transformations extends to a finite group of
free symmetries of $\mathbb S^3$ preserving the 2-bridge knot or link
with exterior $\hat M^3$.

The symmetry group of a hyperbolic 2-bridge knot or link is known by
M. Sakuma's
work ([Sa2]), using Thurston's orbifold Theorem (cf.[BP],[CHK]). In
particular, the orientation preserving symmetry subgroup, acting
freely on the exterior of 2-bridge knot or link, has order at
most two. Moreover, it may have order
two only in the case of a 2-bridge link, since a hyperbolic 2-bridge
knot does not admit
a free symmetry ( [GLM], [Ha], [Sa2]). \hfill$\Box$

\bigskip

\section{The proof of Theorem \ref{JSJ}}

A first observation is that the splitting of $G=\pi_1(M^3)$ that
corresponds to the JSJ-decomposition of the 3-manifold $M^3$ is 2-acylindrical unless one of the Seifert pieces is
homeomorphic to $Q^3$. This follows easily form the
following three facts:

\begin{enumerate}

\item All properly embedded annuli in a Seifert manifold different
from $Q^3$ are
either parallel to the boundary or are vertical meaning that the
intersection with the boundary components are the fibre, this implies
that only the element corresponding to the fibre can fix two edges
emanating at a vertex in the Bass-Serre tree corresponding to a
Seifert piece.

\item For the hyperbolic (acylindrical) pieces there is no such
element at all.

\item The fibre is never glued to the fibre at an essential JSJ-torus
separating
two Seifert pieces.

\end{enumerate}

We now proceed with the proof of Theorem~\ref{JSJ} for the case where
the JSJ-decomposition does not contain any pieces that are
homeomorphic to $Q^3$. We
first consider the case where the JSJ-decomposition of the given
3-manifold has a
separating torus and then where it does not. We conclude by dealing
with the case
that there are pieces of type $Q^3$.

\subsection{The JSJ-decomposition has a separating torus $T$ and no
piece of type
$Q^3$}\label{noq3}

The torus $T$ splits the 3-manifold in two pieces $M_A$ and $M_B$.
The fundamental group $G=\pi_1(M^3)$ then splits as an amalgamated
product $A*_CB$ with amalgam $C\cong \mathbb Z\oplus\mathbb Z$ and
$A$ and $B$ the fundamental groups of $M_A$ and $M_B$. This
amalgamated product satisfies the condition of Lemma~\ref{amalgam}.
We therefore have to investigate the following different cases, where
$\{g,h\}$ is a generating set of $G$. It is clear that we can always
assume $g$ to be primitive.

\medskip\noindent {\bf 1. Case:} $g\in A$ and $g^n\in C$ and $h\in
B-1$.

\smallskip\noindent It is clear that $g\in A-C$ since otherwise
$\langle g,h\rangle \in B$ which implies that $g$ and $h$ do not
generate $G$. In particular it follows that $g$ lies in the
fundamental group of the piece $M_{A_1}$ of the
JSJ-decomposition of $M_A$ that contains $T$ since the splitting is 2-acylindrical.
Since $g\notin C$ this implies that $M_{A_1}$ is Seifert and that $g$
is a root of the fibre. This also implies that $M_A=M_{A_1}$ since
$\langle g,h\rangle\subset \pi_1(M_{A_1}\cup_T M_B)$ and therefore
$M_A$ cannot contain a piece besides $M_{A_1}$. We define $C'=\langle
g^n,h\rangle\cap C\subset B$. Now by Lemma~\ref{top1} we get that
$\langle g,C'\rangle\cap C=C'$. Proposition~\ref{controlled} (1)
therefore implies that the induced splitting of $\langle g,h\rangle $
is $\langle
g,C'\rangle *_{C'}\langle g^n,h\rangle$. Since $g$ and $h$ generate
$G$ this implies that $C'=C$, that $A=\langle g,C\rangle$ and that
$B=\langle g^n,h\rangle$. By Lemma~\ref{top2} $A=\langle g,C\rangle$
implies that $M_A$ is a
Seifert manifold with base $D(p,q)$ or $A(p)$.

It remains to analyze what the manifold $M_B$ can be. If $M_B$ has a
trivial JSJ-decomposition and the only piece is hyperbolic we cannot
say anything more, this
puts us into situation 1 of Theorem~\ref{JSJ}.

If $M_B$ consists of a single Seifert piece then by Lemma~\ref{top5}
$M_B$ must either be the exterior of a 1-bridge knot in a lens space
and $g^n$ corresponds to a meridian of the knot or a Seifert space
with base space $A(p)$ and $g^n$ corresponds to
a curve that has intersection number one with the fibre. This puts us
in situation 2 and 3, respectively, of Theorem~\ref{JSJ}.

If $M_B$ has a non-trivial JSJ-decomposition then we study the action
of $\langle
g^n,h\rangle$ on the Bass-Serre tree of the corresponding splitting
of $B$. It is
clear that $g^n$ acts with a fixed point since it lies in the vertex
group that corresponds to the piece $M_{B_1}$ of the
JSJ-decomposition of $M_B$ that has $T$ as a boundary
component. It is also easy to see that $T_{g^n}$ consists of a single
point: This
is clear if $M_{B_1}$ is acylindrical (hyperbolic). If $M_{B_1}$ is
Seifert it follows from the fact that no power of $g^n$ corresponds
to a power of the fibre of $M_{B_1}$, otherwise $T$ wouldn't be a
torus of the JSJ-decomposition, but merely an essential torus in a
Seifert
piece. By the remark after Theorem~\ref{KW} we can assume that either
$T_{g^n}\cap T_{h}\neq \emptyset$ or that $T_{g^n}\cap h T_{g^n}\neq
\emptyset$. In the second case $h$ must fix the same vertex as $g^n$
and therefore lie in the same vertex group as $g^n$. This implies
that $\langle g^n,h\rangle$ does only intersect the vertex group that
corresponds to $M_{B_1}$ which implies that $g^n$ and $h$ do not
generate $B$ since $M_B$
was assumed to have non-trivial JSJ-decomposition, a contradiction.

It remains to verify the case $T_{g^n}\cap T_{h}\neq \emptyset$:
Since $g^n$ and $h$ act with a fixed point it follows that the
JSJ-decomposition of $M_B$ does not contain a non-separating torus
since otherwise $g^n$ and $h$ cannot generate $B$ by the argument
given in the proof of Lemma~\ref{HNN1}. $T_{g^n}\cap T_{h}\neq
\emptyset$ implies that
a power $h^m$ of $h$ fixes the single vertex $v$ of $T_{g^n}$. By the
above argument we can assume that $h$ does not fix $v$. Since the
action is 2-acylindrical this implies that $h$ fixes a
vertex $w$ that is in distance one from $v$. As before we argue that
the piece $M_{B_2}$ of the JSJ-decomposition of $M_B$ corresponding
to $w$ is a Seifert piece and $h$ is a root of the fibre. Let $T'$ be
the torus of the JSJ-decomposition of $M_B$ that separates $M_{B_1}$
and $M_{B_2}$ and let $C'$ the corresponding edge group. As before we
now see
that the induced splitting of $\langle g^n,h\rangle$ is $\langle
g^n,h^m\rangle *_{C'_1}\langle C'_1,h\rangle$ where
$C'_1=C'\cap\langle g^n,h^m\rangle$ and that
we must have $C'_1=C'$. Now this implies that $M_{B_1}$ is generated
by two element that come from subgroups that correspond to two
different torus boundary components and is therefore a 2-bridge link
by either Lemma~\ref{top3} or Lemma~\ref{hyp1} depending whether the
piece is Seifert or hyperbolic. Furthermore the two generators must
correspond to meridians. As before we further see that $M_{B_2}$ must
be of type $D(p,q)$ or $A(p)$ by Lemma~\ref{top2}.
This puts us in situation 6 of Theorem~\ref{JSJ}.

\medskip\noindent {\bf 2. Case:} $g\in A$ and $g^n\in C$ (we assume
that $n$ is chosen minimal with this property) and $h=ab$ with $a\in
A-C$, $b\in B-C$ and $a^{-1}g^ma\in C$ for some $m\in \mathbb N$.

\smallskip Let $M_{A_1}$ be the piece of the JSJ-decomposition of
$M_A$ that contains
$T$. It is clear that $g^{nm}\in C$ and that $a^{-1}g^{nm}a\in C$.
The properly immersed annulus corresponding to $a^{-1}g^{nm}a=c\in C$
can be homotoped into the piece $M_{A_1}$. This implies that $\langle
g,h\rangle\subset \pi_1(M_{A_1}\cup_T M_B)$ and therefore
$M_A=M_{A_1}$. Now $M_A$ cannot be acylindrical since otherwise the
annulus must be boundary parallel which implies that $a\in C$ which
contradicts our assumptions. Thus $M_A$ is Seifert. Since we assume
that $M_A$ is not of type $Q^3$ this implies that $g^{mn}$ is a power
of the fibre of $M_A$. After replacing $g$ by a generator of the
cyclic group $\langle
g,f_A\rangle$ we can assume that g is a root of the fibre $f_A$ of
$M_A$. In particular we have $f_A=g^n$ and $a^{-1}g^na=g^{\pm n}$ and
therefore $h^{-1}g^nh=b^{-1}g^{\pm n}b$.

\smallskip\noindent We first show that we can restrict ourselves to
the case that $\langle g^n,b^{-1}g^{n}b\rangle$ is not free. Suppose
that $\langle g^n,b^{-1}g^{n}b\rangle$ is free. We look at the
Bass-Serre tree with
respect to the decomposition $G=A*_CB$. Let $w$ be the vertex fixed
under the action of $A$ and $v$ be the vertex fixed under the action
of $B$. It is clear that the vertex $z=h^{-1}w=b^{-1}a^{-1}w=b^{-1}w$
is different from $w$ and in distance $1$ from $v$. We denote the
edge $[v,w]$ by $e_1$ and the edge $[v,z]$ by $e_2$. It is
clear that that $b^{-1}e_1=e_2$. We define $G_w=\langle g\rangle$,
$G_z=\langle h^{-1}gh\rangle$ and $G_v=\langle
g^n,b^{-1}g^nb\rangle$. The freeness of $G_v$ and the minimality of
$n$ guarantees
that $G_v\cap \Stab e_1=\langle g^n\rangle=G_w\cap \Stab e_1$, that
$G_v\cap \Stab e_2=\langle b^{-1}g^nb\rangle=G_z\cap \Stab e_2$ and
that $e_1$ and $e_2$ are not $G_v$-equivalent. Since $hv\neq v$ we
have that $he_2=abb^{-1}e_1=ae_1\neq
e_1$. Now $ae_1$ and $e_1$ are not $\langle g\rangle$-equivalent
unless $ac\in\langle g\rangle$ for some $c\in C$ in which case we can
replace $h$ by an
element of $B$ after left multiplication with a power of $g$. This
however puts us into the first case, we can therefore assume that
$ae_1$ and $e_1$ are not $\langle g\rangle$-equivalent. It follows
that all the conditions of Proposition~\ref{controlled} (2) are
fulfilled, i.e. the induced splitting of $\langle g,h\rangle$ has two
edges and vertices with cyclic edge groups. This however means that
$\langle g,h\rangle\neq
G$.

We next show that the JSJ-decomposition of $M_B$ must be trivial,
i.e. that $M_B=M_{B_1}$ is geometric, where $M_{B_1}$ is the piece of
the JSJ of $B$ containing $T$. It clearly suffices to show that the
subgroup $\langle C,b\rangle$ lies in the subgroup $B_1$ of $B$
corresponding to $M_{B_1}$. We study the action of $B$ on the
Bass-Serre tree associated to the splitting of $B$
that corresponds to the JSJ-decomposition of $M_B$. Clearly $g^n$
acts with a fixed
point and since $g^n$ is not the fibre of a Seifert piece of $M_B$ it
follows that it is
not conjugate to an element of one of the edge groups of the
decomposition of $B$, nor
is one of its powers. It follows that $T_{g^n}$ consists of a single
vertex and so does $T_{h^{-1}g^nh}=h^{-1}T_{g^n}$. Now by Lemma~2.1
of [KW] either $\langle
g^n,h^{-1}g^nh\rangle$ is free or $T_{g^n}\cap T_{h^{-1}g^nh}\neq
\emptyset$. Since we already dealt with the first case we can assume
that the second case holds. It follows that $h$ fixes the vertex
which is fixed under the action of $g^n$, hence $b\in B_1$ which
proves our assertion.

\smallskip\noindent We distinguish the cases that $M_B$ is a Seifert
manifold and that $M_B$ is hyperbolic.

\smallskip\noindent Suppose that $M_B$ is a Seifert manifold. We have
to investigate the following cases because of Lemma~\ref{top6}.

\smallskip\noindent (1) $M_B$ is the complement of a 2-bridge knot,
$g^n$ corresponds to a meridian, $a^{-1}g^na=g^{\pm n}$, i.e.
$g^n=f_A^{\pm 1}$ and $\langle g^n,h^{-1}g^nh\rangle=B$.

\noindent We get that $\langle g,h=ab\rangle=\langle g,a,b\rangle$
since $b\in B=\langle g^n,h^{-1}g^nh\rangle$. It follows that
$\langle g,h\rangle=G$ ifandonlyif $A=\langle g,a,C\rangle$ where $g$
is a
root of the fibre (or the fibre itself).
This clearly implies that $M_A$ is a Seifert space with base space
$M\ddot o(p)$, $M\ddot o(p,q)$, $D(p,q)$, $D(p,q,r)$, $A(p)$,
$A(p,q)$, $\Sigma$,
$\Sigma(p)$ or the once punctured M\"obius band with at most one cone
point which
puts us into situation~4 of Theorem~\ref{JSJ}. This follows since the
group obtained by killing the boundary curve must be two generated
with one elliptic generator. The only two-generated closed 2-orbifold
group not in this list is the
group of the base $D(2,2,2,2l+1)$ which however has no generating
pair that contains an elliptic element \cite{PRZ}.

\smallskip\noindent (2) $M_B$ is a Seifert manifold with base
$D(2,2l+1)$, $g^n$
corresponds to a curve that has intersection number one with the
fibre of $M_B$,
$a^{-1}g^na=g^{\pm n}$ and $\langle g^n,h^{-1}g^nh\rangle$ maps
surjectively onto
the base group and $\langle g^n,h^{-1}g^nh\rangle\neq B$.

\noindent Note that $|B:\langle g^n,h^{-1}g^nh\rangle|=k<\infty$
where $k>1$ is such that $\langle f_B\rangle\cap \langle
g^n,h^{-1}g^nh\rangle=\langle f_B^k\rangle$. Let $m>2$ be a prime
such that $m|k$. We add the generator $f_B^m$ to the generating pair
$\{g,h\}$. We clearly have $\langle f_B\rangle\cap \langle
g^n,h^{-1}g^nh,f_B^m\rangle=\langle f_B^m\rangle$. Since $\langle
g^n,h^{-1}g^nh\rangle$ maps surjectively onto the base group we know
that $f_B^lb\in \langle g^n,h^{-1}g^nh\rangle$ for some $l\in\mathbb
Z$. We then rewrite the reduced form of $h=ab$ in the form
$h=a(f_B^{-l}f_B^l)b=(af_B^{-l})(f_B^lb)=a'b'$ and denote $a'$ and
$b'$ again by
$a$ and $b$. We then have that $\langle g,h=ab,f_B^m\rangle=\langle
g,a,b,f_B^m\rangle$. Let
now $B'=\langle a,g,f_B^m\rangle$. Now either $B'\cap C=\langle
g^n,f_B^m\rangle$
or $B'\cap C=C$ since there is no intermediate group between $\langle
g^n,f_B^m\rangle$ and $C$ since $m$ is prime. We now apply
Proposition~\ref{controlled} (1). In the first case we get that the
induced splitting of $\langle g,h,f_B^m\rangle$ is $\langle
g,a,f_B^m\rangle*_{\langle f_B^m,g^n\rangle} \langle
g^n,h^{-1}g^nh,f_B^m\rangle$ and
in particular that $\langle g,h,f_B^m\rangle\cap C=\langle
g^n,f_B^m\rangle\neq C$
which contradicts the fact that $G=\langle g,h,f_B^m\rangle$. In the
second case we get
that the induced splitting is of type $\langle g,a,f_B^m\rangle *_C
B$. If $g$ and $h$ generate $G$ this means that $A=\langle
g,a,f_B^m\rangle$. This implies that the quotient of the base group
by the $m$th power of the boundary curve is two generated where one
generator is elliptic. This implies that the original base was of one
of the types $M\ddot o$, $M\ddot o(p)$, $D(p,q)$ or $A(p)$. The first
two cases yield situation~5 of Theorem~\ref{JSJ}, the last two cases
yield situation~2. Note that these groups are in fact $2$-generated
by \cite{Ko2} and \cite{W2}.

\smallskip\noindent (3) $M_B$ is a Seifert manifold with base
$D(2,2l)$, $g^n$ corresponds to a curve that has intersection number
one with the fibre of $M_B$. It follows from Lemma~\ref{top6} that
$g^n$ and $h^{-1}g^nh=b^{-1}g^nb$ are not conjugate in $\langle
g^n,h^{-1}g^nh\rangle$, nor are their images under $\pi$ conjugate in
the image of $\langle g^n,h^{-1}g^nh\rangle$. Hence they are also not
conjugate in $\langle g^n,h^{-1}g^nh,f_B\rangle$ when we add the
generator $f_B$. In particular $b\notin \langle
g^n,h^{-1}g^nh,f_B\rangle$. We distinguish the cases where $a\in
\langle
g,f_B,af_Ba^{-1}\rangle$ and $a\notin\langle
g,f_B,af_Ba^{-1}\rangle$.

In the first case we have that $G=\langle g,h=ab,f_B\rangle=\langle
g,a,f_B,B\rangle$ and Proposition~\ref{controlled} (1) yields that
$G=\langle g,h,f_B\rangle=\langle g,a,f_B\rangle *_{\langle
g^n,f_B\rangle}\langle g^n,f_B,b\rangle=\langle g,a,C\rangle
*_C\langle C,b\rangle$. In particular $\langle g,a,C\rangle=\langle
g,f_B,af_Ba^{-1}\rangle=A$ and $\langle C,b\rangle=B$. Let $N_A(C)$
be the normal closure of $C$ in $A$. Now $A/N_A(C)$ must be generated
by one element ($f_B$ and $af_Ba^{-1}$ lie in the kernel), namely the
image of $g$. This however implies that the base is of type $D(p,q)$
or $A(p)$ since the only other cyclic 2-orbifold group is the
projective plane with at most one singularity, but here the generator
is not elliptic. This puts us into situation~2 of Theorem~\ref{JSJ}.
Again we know that $G$ is in fact $2$-generated since the manifold is
of genus $2$ by \cite{Ko2}.

In the second case we choose $v$, $w$, $z$, $e_1$ and $e_2$ as in the
beginning of case (2) above. We define $G_v=\langle
g,f_B,af_Ba^{-1}\rangle$, $G_w=\langle
f_B,g^n,b^{-1}g^nb\rangle $ and $G_z=\langle
b^{-1}g^nb,b^{-1}a^{-1}f_Bab,b^{-1}a^{-1}gab\rangle$. It is clear
that $e_2=be_1$
and that $he_2=ae_1$. The facts that $b\notin \langle
f_B,g^n,b^{-1}g^nb\rangle$
and that $a\notin \langle g,f_B,af_Ba^{-1}\rangle$ imply that $e_1$
and $e_2$ are
$G_w$-inequivalent and that $e_1$ and $he_2$ are $G_v$-inequivalent.
Proposition~\ref{controlled} (4) then implies that the induced
splitting of $\langle g,h,f_B\rangle$ has two edge groups and two
vertex groups, which implies
that it is not the induced splitting of $G$, i.e. $\langle
g,h,f_B\rangle\neq G$
and therefore $\langle g,h\rangle \neq G$.

\medskip\noindent Suppose now that $M_B$ is hyperbolic.

\medskip\noindent Since $U:=\langle
g^n,h^{-1}g^nh\rangle=\langle g^n,b^{-1}g^nb\rangle$ is neither free
nor Abelian it follows from Proposition~\ref{twobridge} that either
$M_B$ is the complement of a 2-bridge knot or link, that $U=B$ and
that $g^n$ and $b^{-1}g^nb$ correspond to meridians or that
$U$ is a subgroup of $B$ of index $2$ and that the covering space
$\hat M_B$ corresponding to $U$ is homeomorphic to the exterior of a
2-bridge link in $\mathbb S^3$ (with two components).

\medskip In the first case we argue precisely as in case (1) above.
In the second
case we can assume that the
covering $p:\hat M_B\to M_B$ is a homeomorphism when restricted to a
boundary component since the degree of the covering is $2$ and since
the two boundary components of $\hat M^3$ get mapped onto the same
boundary component of $M^3$. It follows
that $C\cap U=C$, $h^{-1}Ch\cap U=h^{-1}Ch$ and $b\notin U$ since the
two meridians are not conjugate in $U$. It follows that $B=\langle
g^n,b\rangle$ where $g^n$ is parabolic. We can now argue precisely as
in case (3) above and either obtain that $\langle g,h\rangle\neq G$
or that $M_A$ is a Seifert space over $D(p,q)$ or $A(p)$. Since the
fibre $f_{M_A}=g^n$of $M_A$ gets identified with the parabolic
generator of $B$ this puts us into situation~1 of Theorem~\ref{JSJ}.

\medskip\noindent {\bf 3. Case:} $g\in A$ and $g^n\in C$ and
$h=bab^{-1}$ with $a^m\in C-1$ for some $m\in \mathbb N$.

\smallskip\noindent We can assume that $g\notin C$ since we are
otherwise in the
first case after conjugation with $b$. As in the case before it
follows that $M_A$ is Seifert and that we can assume that $a,g\in
A-C$ are roots of the fibre $f_A$.

We distinguish the cases that $bc\in
\langle g^n,ba^mb^{-1}\rangle\subset B$ for some $c\in C$ and that
$bc\notin \langle g^n,ba^mb^{-1}\rangle$ for all $c\in C$.

\smallskip\noindent (a) Suppose that $bc\in \langle
g^n,ba^mb^{-1}\rangle\subset B$ for some $c\in C$. After rewriting
the reduced form of $h=bab^{-1}$ as $(bc)(c^{-1}ac)(c^{-1}b^{-1})$ we
have that $b\in \langle g^n,ba^mb^{-1}\rangle$.
It follows that $\langle g,h\rangle=\langle g,b,a\rangle$. Since $g$
and $a$ are both roots of the fibre it follows that $g^n=a^m=f_A$.
This implies that $ba^mb^{-1}=(ba)g^n(a^{-1}b^{-1})$. We show that
$\langle g,h\rangle=\langle g,ab\rangle$ which puts us into the
second case. Since
$ba^mb^{-1}=bf_Ab^{-1}=(ba)f_A(ba)^{-1}=(ba)g^n(ba)^{-1}\in \langle
g,ab\rangle$ it follows that $b\in\langle g,ab\rangle$ and therefore
also $a\in \langle g,ab\rangle$. We have shown that $\langle
g,ab\rangle=\langle a,b,g\rangle=\langle g,h\rangle$.

\smallskip\noindent (b) Suppose that $bc\notin \langle
g^n,ba^mb^{-1}\rangle$ for all $c\in C$. We study the action on the
Bass-Serre tree associated to
the splitting $G=A*_CB$. Let $x$ be the vertex fixed under the action
of $A$, $y$
be the vertex fixed by $B$ and $z=bx$ be the vertex fixed under the
action of $bab^{-1}$. We further define $e_1=[x,y]$ and $e_2=[y,z]$.
Note that $\Stab e_1=C$ and that $\Stab e_2=bCb^{-1}$. We define
$G_y=\langle g^n,ba^mb^{-1}\rangle$, $G_x=\langle g,C\cap G_y\rangle$
and $G_z=\langle h,bCb^{-1}\cap G_z\rangle$. Lemma~\ref{top1}
guarantees that $G_x\cap C=G_y\cap C$ and that $G_y\cap
bCb^{-1}=G_z\cap bCb^{-1}$. We further have that $e_1$ and $e_2$ are
$G_y$-inequivalent since $bc \notin G_y$ for all $c\in C$. This
implies
that all hypothesis of Proposition~\ref{controlled} (2) are
fulfilled, i.e. that
the induced splitting of $\langle g,h\rangle$ has two edge groups.
This however implies that $\langle g,h\rangle\neq G$.

\subsection{The JSJ-decomposition has a non-separating torus}

We can assume that there is no separating torus because there is no
3-manifold that has 2-generated fundamental group such that the JSJ
contains a separating and a
non-separating torus. In the case when there is no piece of type
$Q^3$ this follows from section~\ref{noq3}, in the case with a piece
of type $Q^3$ we will see this in section~\ref{q3}. This implies that
the
JSJ-graph is homeomorphic to a circle and we only have to look at the
cases where
we have one or two pieces because of Lemma~\ref{HNN1}.

\medskip\noindent {\bf 1. case: } The JSJ-decomposition has one piece
$N$ with one non-separating torus.

\smallskip\noindent Let $A$ be the fundamental group of $N$. We can
then write $G$ as the HNN-extension $\langle
A,t|tc_1t^{-1}=\phi(c_1)\hbox{ for }c_1\in C_1\rangle$ where
$\phi:C_1\to C_2$ is the isomorphism between the two torus subgroups
induced by the JSJ-decomposition of the manifold. Because of
Lemma~\ref{HNN3} we can
assume that there exists a generating pair $\{g,h\}$ such that $g\in
A$ and that
$g^n\in C_1$ and that $h=at$ for some $a\in A$. Note that
$hg^nh^{-1}\in A$. We first look at the case that $N$ is Seifert and
then at the case that $N$ is a hyperbolic piece. In the first case we
distinguish the cases that neither $g^n$ nor $hg^nh^{-1}$ corresponds
to the fibre of $N$ and that either $g^n$ or
$hg^nh^{-1}$ corresponds to the fibre of $N$. We can clearly choose
$g$ such that $g$ generates a maximal cyclic subgroup and $n$ minimal
such that $g^n\in C_1$.

\smallskip\noindent (a) Suppose that $N$ is Seifert and neither $g^n$
nor $hg^nh^{-1}$ correspond to the fibre. This implies that $g\in
C_1$, otherwise $g$
would be a root of the fibre and $g^n=f_N$. If $\langle
g,hgh^{-1}\rangle$ is free then $\langle g,hgh^{-1}\rangle\cap
C_1=\langle g\rangle$ and $h^{-1}\langle g,hgh^{-1}\rangle h\cap
C_1=\langle h^{-1}gh,g\rangle\cap C_1=\langle g\rangle$.
Proposition~\ref{controlled} (3) therefore implies that the induced
splitting of $\langle g,hgh^{-1}\rangle$ is an
HNN-extension with cyclic edge group which implies that $\langle
g,hgh^{-1}\rangle\neq G$. It therefore follows from Lemma~\ref{top3}
that the base $\B$ of $N$ is of type $A(p)$. If $\langle
g,hgh^{-1}\rangle$ is not free we can by Lemma~\ref{top3} assume that
$\langle g,hgh^{-1}\rangle$ maps surjectively onto $F$ and that $g$
and $hgh^{-1}$ correspond to curves on the boundary that have
intersection number $1$
with the fibre.

If $\langle g,hgh^{-1}\rangle =A$ then $N$ is a 2-bridge link by
Lemma~\ref{top3} and $g$ and $hgh^{-1}$ correspond to the meridians,
in particular $M^3$ is obtained from $N$ by identifying the boundary
components such that the meridians get identified, i.e. we are in
situation~8 of Theorem~\ref{JSJ}.

If $\langle g,hgh^{-1}\rangle \neq A$ then $|A:\langle
g,hgh^{-1}\rangle|=k\le\infty$ where $k$ is such that $\langle
f\rangle\cap\langle g,hgh^{-1}\rangle =\langle f^k\rangle$. This is
clear since $\langle g,hgh^{-1}\rangle$ maps surjectively onto the
base group. It follows that $C_1\cap\langle g,hgh^{-1}\rangle=\langle
f^k,g\rangle$ and that $C_1\cap h^{-1}\langle g,hgh^{-1}\rangle
h=C_1\cap \langle h^{-1}gh,g\rangle =\langle f^k,g\rangle$. It then
follows from Proposition~\ref{controlled} (3) that the induced
splitting of $\langle g,h\rangle$ has one edge group which is a
proper subgroup of the edge group of the splitting of $G$ which
implies that $G\neq\langle g,h\rangle$.

\smallskip\noindent (b) $N$ is Seifert and either $g^n$ or
$hg^nh^{-1}$ correspond to the fibre. Without loss of generality we
can assume that $g^n$ corresponds to a fibre, i.e. $g$ maps onto an
elliptic element of the base group
and $g^n=f$. Note that $h^{-1}fh=t^{-1}a^{-1}fat=t^{-1}f^{\pm 1}t\in
t^{-1}C_2t=
C_1$. In particular we have $\langle f,h^{-1}fh\rangle\subset C_1$.
We distinguish the cases where $C_1=\langle f,h^{-1}fh\rangle$ and
$C_1\neq\langle f,h^{-1}fh\rangle$.

Suppose that $C_1=\langle f,h^{-1}fh\rangle$. Note that $C_1\subset
\langle g,hfh^{-1},h^{-1}fh\rangle\subset A$ and $C_1\subset
h^{-1}\langle g,hfh^{-1},h^{-1}fh\rangle h=\langle
h^{-1}gh,f,h^{-2}fh^2\rangle\subset h^{-1}Ah$.
Proposition~\ref{controlled} (3) therefore implies that the induced
splitting of $\langle g,h\rangle$ is a HNN-extension with base
$\langle g,hfh^{-1},h^{-1}fh\rangle$ and edge group $C_1$. Note that
$C_1=\langle f,h^{-1}fh\rangle$ implies that $hfh^{-1}$ and
$h^{-1}fh$ correspond to boundary
curves that have intersection number one with the fibre. It follows
that $G=\langle g,h\rangle$ if and only if $F$ is generated by one
elliptic element and two
elements that correspond to boundary curves. The only such orbifolds
are $A(p)$, $A(p,q)$, $\Sigma$ and $\Sigma(p)$, i.e. $M^3$ is
obtained from a Seifert manifold with base $A(p)$, $A(p,q)$, $\Sigma$
or $\Sigma(p)$ where the boundary components are glued such that the
fibre on one component is glued with a curve on the other component
that has intersection number one with the fibre $f$. This puts us
into situation~9 of
Theorem~\ref{JSJ}.

Suppose that $C_1\neq\langle f,h^{-1}fh\rangle$. This means that
$C_1\cap \langle
f,h^{-1}fh\rangle=\langle f,m^k\rangle$ where $m$ corresponds to a
curve with intersection number one with the fibre and $k\ge 2$. In
particular we have $\langle f,m^k\rangle\subset\langle
g,hfh^{-1},h^{-1}fh\rangle$ and $\langle f,m^k\rangle\subset
h^{-1}\langle g,hfh^{-1},h^{-1}fh\rangle h=\langle
h^{-1}gh,f,h^{-2}fh^2\rangle$. It suffices to show that $C_1\cap
\langle g,hfh^{-1},h^{-1}fh\rangle=C_1\cap h^{-1}\langle
g,hfh^{-1},h^{-1}fh\rangle h=\langle f,m^k\rangle$ since we then see
as in the case before that the induced
splitting of $\langle g,h\rangle$ is a HNN-extension with edge group
$\langle f,m^k\rangle$ which implies that $\langle g,h\rangle\neq G$.
We show that $C_1\cap \langle g,hfh^{-1},h^{-1}fh\rangle=\langle
f,m^k\rangle$, the argument for the other statement is analogous.
This can be seen by quotienting the normal
closure $\bar N$ of $\langle g,hfh^{-1},h^{-1}fh\rangle$ out of $A$.
The quotient map maps $A$ onto the fundamental group of the orbifold
that is obtained from the base
space by replacing the two boundary components by cone points of
order $k$. It is clear that $m$ gets mapped onto an elliptic element
of order $k$. This implies
that $\langle m\rangle \cap \bar N=\langle m^k\rangle$ and therefore
$\langle m\rangle \cap \langle g,hfh^{-1},h^{-1}fh\rangle= \langle
m^k\rangle$ which proves the assertion.

\smallskip\noindent (c) Suppose that $N$ is hyperbolic. Since $C_1$
is malnormal
in $A$ this implies that $g\in C_1$ since $g^n\in C_1$. As in case
(a) we see that $\langle
g,hgh^{-1}\rangle$ cannot be free. Since $\langle g,hgh^{-1}\rangle$
is also not Abelian it follows from Proposition~\ref{twobridge} that
either $N$ is the complement of a 2-bridge link and $\langle
g,hgh^{-1}\rangle=A$ or that $\langle g,hgh^{-1}\rangle$ is a
subgroup of $A$ of index two and the corresponding covering is the
exterior of a 2-bridge link. In the first case $g$ and $hgh^{-1}$
correspond to meridians by Proposition~\ref{twobridge} which puts us
into situation~8 of Theorem~\ref{JSJ}. In the second case we can
argue as in the
second part of (a) to show that $\langle g,h\rangle\neq G$. This is
possible since
the degree of the covering is $2$ when restricted to either boundary
components.

\medskip
\noindent {\bf 2. case: } The JSJ-decomposition has two pieces and
two non-separating tori. Let $M_A$ and $M_B$ be the two pieces of the
decomposition and denote $\pi_1(M_A)$ by $A$ and $\pi_1(M_B)$ by $B$.
Let $C_1$ and $C_2$ be the subgroups of $A$ corresponding to the
boundary components of $M_A$ glued with
boundary components of $M_B$. Let $C_3$ and $C_4$ be the
corresponding subgroups of $B$.
Suppose further that $\phi_1:C_1\to C_3$ and $\phi_2:C_2\to C_4$ are
the isomorphisms that are induced by the gluing. Then
$G=(A*_{C_1=C_3}B)*_{C_2=C_4}$,
i.e. $G=\langle A,B,t|\phi_1(c_1)=c_1,\phi_2(c_2)=tc_2t^{-1}\rangle$.
Since this
splitting is 2-acylindrical Lemma~\ref{HNN2} guarantees
(possibly after exchanging $A$ and $B$) that $g\in B$, $g^n\in C_3-1$
and $g^m\in bC_4b^{-1}-1$ for some $b\in B$ and $n,m\in \mathbb N$
and $h=bta$ for some $a\in
A$.

It is clear that $g^{nm}\in C_3-1$ and $g^{mn}\in bC_4b^{-1}$. It
follows from the annulus Theorem that
there exists an essential annulus in $M_B$. This implies that $M_B$
is Seifert and
that $g^{nm}$ corresponds to a power of the fibre. After replacing
$g$ with a generator of a maximal cyclic subgroup containing the
original $g$ we can therefore assume that $n=m$ and $g^n=f_B$.

We now look at the action of $G$ on
the Bass-Serre tree corresponding to the above splitting. Let $w$ be
the vertex fixed by $B$, $v$ be the vertex fixed by $A$ and
$z=h^{-1}w=a^{-1}t^{-1}b^{-1}w=a^{-1}t^{-1}w$ be the vertex fixed by
$(bta)^{-1}B(bta)=(ta)^{-1}B(ta)$. Note that $e_1=[v,w]$ and
$e_2=[v,z]$ are $A$-inequivalent and that $e_1$ and $(bta)e_2$ are
$B$-inequivalent since they project onto distinct edges of the
quotient graph. This implies that for three groups $G_w$, $G_v$ and
$G_z$ fixing $w$, $v$ and $z$, respectively, we only need
to show that $G_v\cap \Stab e_1=G_w\cap\Stab e_1$, $G_v\cap \Stab
e_2=G_z\cap\Stab e_2$ and that $hG_wh^{-1}=G_z$ in order to verify
that the hypothesis of Proposition~\ref{controlled} (4) are
fulfilled.

Note that $\langle
f_B,h^{-1}f_Bh\rangle=\langle f_B,(bta)^{-1}f_B(bta)\rangle=\langle
f_B,a^{-1}t^{-1}f_B^{\pm 1}ta\rangle \subset A$ and that $f_B$ and
$h^{-1}f_Bh$ correspond to elements of different boundary components
of
$M_A$. Using the analysis of (a) and (c) above we have to distinguish
three cases (i) where $\langle f_B,h^{-1}f_Bh\rangle$ is free, (ii)
where $\langle
f_B,h^{-1}f_Bh\rangle=A$ and $M_A$ is a 2-bridge link and $f_B$ and
$h^{-1}f_Bh$ correspond to meridians and (iii) where $\langle
f_B,h^{-1}f_Bh\rangle$ is of finite index in $A$ and
$|C_1:C_1'|=|C_4:C_4'|=|A:\langle f_B,hf_Bh^{-1}\rangle|$ where
$C_1'=C_1\cap \langle f_B,hf_Bh^{-1}\rangle$ and $C_4'=C_4\cap
a\langle f_B,hf_Bh^{-1}\rangle a^{-1}$.

\smallskip\noindent (i) If $A'=\langle f_B,h^{-1}f_Bh\rangle\subset
A$ is free then clearly $A'\cap C_1=\langle f_B\rangle$ and $A'\cap
h^{-1}C_2h=\langle h^{-1}f_Bh\rangle$. This implies that if we define
$G_w=\langle g\rangle$, $G_v=\langle f_B,h^{-1}f_Bh\rangle$ and
$G_z=\langle h^{-1}gh\rangle$ the hypothesis of
Proposition~\ref{controlled} (4) are fulfilled and therefore the
induced splitting of $\langle g,h\rangle$ has cyclic edge groups
which implies that $\langle g,h\rangle\neq G$.

\smallskip\noindent (ii) If $\langle f_B,h^{-1}f_Bh\rangle=A$ then in
particular
$C_1\subset \langle f_B,h^{-1}fh\rangle$ and $a^{-1}C_4a\subset
\langle f_B,h^{-1}f_Bh\rangle$. We define $G_w=\langle
g,C_1,bC_2b^{-1}\rangle$, $G_v=\langle f_B,h^{-1}f_Bh\rangle=A$ and
$G_z=h^{-1}G_wh$. This implies that the
induced splitting of $\langle g,h\rangle$ consist of two vertex
groups where one
equals $A$ and the other to $G_w$ and the edge groups are the edge
groups of the
original group. This means that $\langle g,h\rangle=G$ ifandonlyif
$G_w=B$.
This however
implies that $B$ is generated by the subgroups corresponding to the
boundary components and a root of the fibre, i.e. the base is
generated by two boundary curves and one elliptic element, as before
we see that this implies that the base
space is of type $A(p)$, $A(p,q)$, $\Sigma$ or $\Sigma(p)$. This
means that $M^3$ is
obtained from a Seifert manifold with base space $A(p)$, $A(p,q)$,
$\Sigma$ or $\Sigma(p)$ and a 2-bridge link complement where the
boundary components are glued such that the fibre is glued to
meridians of the 2-bridge link complement.
This puts us into situation~7 of Theorem~\ref{JSJ}.

\smallskip\noindent (iii) Let $C_2'=tC_4't^{-1}$. We define
$G_w=\langle g,C_1',bC_2'b^{-1}\rangle$, $G_v=\langle
f_B,hf_Bh^{-1}\rangle$ and $G_z=h^{-1}\langle
g,C_1',bC_2'b^{-1}\rangle h=\langle
h^{-1}gh,hC_1h^{-1},a^{-1}C_4'a\rangle$. The same argument as in case
1 (c) shows
that $G_w\cap\Stab e_1=C_1'$ and that $G_z\cap\Stab e_2=a^{-1}C_4'a$,
i.e. that the hypothesis of Proposition~\ref{controlled} (4) is
fulfilled. This implies that one of the vertex groups of the induced
splitting of $\langle g,h\rangle$ is a proper subgroup of $A$ and the
other a proper subgroup of $B$. It follows that $\langle
g,h\rangle\neq G$.

\subsection{Dealing with the existence of pieces of type
$Q^3$}\label{q3}

We start by the following lemma that will be useful in the course of
our investigation.

\begin{Lemma}\label{Q^3&S}
Let $M^3$ be a 3-manifold whose JSJ-decomposition consists of a piece
homeomorphic
to $Q^3$ and a Seifert piece $M_B$ with base $D(p,q)$ or $A(p)$.
Suppose further that
$\pi_1(M^3)=\langle g,h\rangle$ where $g$ is a root of the fibre
$f_B$ of $M_B$. Then the
intersection number of $f_Q'$ and $f_B$ is one.
\end{Lemma}

\noindent{\em Proof} We use the notation as in the beginning of
section~\ref{toplem}. Let $G=A*_C B$ be the splitting of $G$ that
corresponds to the
JSJ-decomposition of $M^3$. Suppose that the intersection number of
$f_Q'$ and $f_B$ is greater than one, i.e. that $f_B=(xy)^n{f_Q'}^m$
for some $n\ge 2$. It follows that
$xf_Bx^{-1}=(xy)^{-n}{f_Q'}^{m+2n}$. We clearly have $\langle
f_B,xf_Bx^{-1}\rangle\subset \langle (xy)^n,f_Q'\rangle=:U\neq C$ and
$U$ is normal in $A$ with $A/U\cong D_{2n}$ and $C/U\cong \mathbb
Z_n$. In particular $g^2=1$ for all $g\in A/U-C/U$.

Since $(xy)^n$ is a $n$th power and since $f_Q'$ has intersection
number $n$ with $f_B$ it follows that the intersection number of any
element of $U$ with $f_B$ is a multiple of $n$. It follows that
$B/N_B(U)$ is the fundamental group of the orbifold obtained by
replacing a boundary curve of the base space of $M_B$ with a cone
point of order $n$. If the base space was $D(p,q)$ the resulting
orbifold is $S^2(p,q,n)$ if it was $A(p)$ it becomes $D(p,n)$.

It follows that $G'=G/N_G(U)$ can again be written as an amalgamated
product $G'=A'*_{C'}B'$ where $A'=A/U$, $C'=C/U$ and $B'=B/N_B(U)$.
Denote the quotient map $G\to G'$ by $\phi$. Since $G'$ is a proper
amalgamated product and therefore not cyclic it follows that
$\phi(g)$ cannot be trivial. It is clear that no non-trivial power of
$\phi(g)$ is conjugate to an element of $C'$ since in a triangle
groups elliptic elements that correspond to different cone points in
the orbifold are not conjugate and the amalgam corresponds to the new
cone points. This implies that $T_{\phi(g)}$ consists of a single
vertex if we look at the action of $G'$ with respect to the splitting
$G'=A'*_{C'}B'$.

Note that $G'$ is not the free product of two cyclic groups. By the
main result of \cite{KW} we can now choose $h$ such that either
$T_{\phi(g)}\cap T_{\phi(h)}\neq \emptyset$ or that $T_{\phi(g)}\cap
hT_{\phi(g)}\neq\emptyset$. The case $T_{\phi(g)}\cap
{\phi(h)}T_{\phi(g)}\neq\emptyset$ cannot occur since then
${\phi(h)}$ would also fix the single vertex of $T_{\phi(g)}$ which
implies that $\langle {\phi(g)},{\phi(h)}\rangle \subset B$.

It follows that $T_{\phi(g)}\cap T_{\phi(h)}\neq \emptyset$. In
particular $\phi(h)$ acts with a fixed point. If $\phi(h)$ is
conjugate to an element of $B'$ then $\langle \phi(g),\phi(h)\rangle$
lies in the kernel of $G'\to G'/B'\cong \mathbb Z_2$, i.e. $G'\neq
\langle \phi(g),\phi(h)\rangle$. We can therefore assume that
$\phi(h)$ is conjugate to an element of $A-C$. Clearly these elements
are of order two and fix no edge. This implies that $T_{\phi(h)}$
consists of a single vertex. It then follows from $T_{\phi(g)}\cap
T_{\phi(h)}\neq \emptyset$ that $T_{\phi(g)}= T_{\phi(h)}$ since
$T_{\phi(g)}$ consists of a single vertex also. This however gives a
contradiction since the vertices correspond to different factors of
the amalgamated product.\hfill$\Box$

\medskip\noindent The preceding lemma allows us to deal with the case
that the JSJ-decomposition has two pieces that are both homeomorphic
to $Q^3$; the following Lemma implies that manifolds of this type
that have 2-generated fundamental group fall into situation~2 of
Theorem~\ref{JSJ}.

\begin{Lemma}\label{2Q^3} Let $M^3$ be a 3-manifold whose
JSJ-decomposition consists of two pieces that are homeomorphic to
$Q^3$. Then $\rank\pi_1(M^3)=2$ ifandonlyif the fibre of one piece has
intersection number one with the fibre of the other piece (when we
look at the pieces as the Seifert space over $D(2,2)$). \end{Lemma}

\noindent{\em Proof} Let $G=A*_C B$ be the splitting of $G$ that
corresponds to the
JSJ-decomposition of $M^3$. Suppose that $G$ is $2$-generated. It
follows from Theorem~\ref{KW} that $G$ is generated by a pair of
elements $\{g,h\}$ such that $g$ acts with a fixed point. It is clear
that $G$ is not conjugate to an element of $C$, since it would
otherwise lie in the kernel of the quotient map $G\to G/C\cong\mathbb
Z_2*\mathbb Z_2$ which clearly implies that $\{g,h\}$ do not generate
$G$. It follows that after conjuration either $g\in A-C$ or $g\in
B-C$, in particular after replacing $g$ with a primitive element we
have that $g^2$ is the fibre of one of the pieces. The assertion now
follows from Lemma~\ref{Q^3&S}.~\hfill$\Box$

\medskip\noindent Next we show that basically the same result as in
Lemma~\ref{amalgam} holds if one of the pieces is homeomorphic to
$Q^3$.

\begin{Lemma}\label{specialcase} Suppose that $M^3$ is a 3-manifold,
$T$ a separating torus of the JSJ-decomposition such that
$M^3=Q^3\cup_T M_B$ where
$M_B$ is not homeomorphic to $Q^3$. Let $G=A*_CB$ be the
corresponding decomposition of $G=\pi_1(M^3)$. Suppose that $\rank
G=2$. Then either $M^3$ falls into case 2 of Theorem~\ref{JSJ} or
there exists a generating
pair $\{g,h\}$ such that either 1 or 2 hold:

\begin{enumerate}
\item $g\in A-C$, $g^2\in C$ and one of the following holds:

(i) $h\in B-C$.

(ii) $h=ab$ with $a\in A-C$, $b\in B-C$.

(iii) $h=bab^{-1}$ with $a\in A-C$ and $a^2\in C-1$.

\item $M_B$ is Seifert, $g\in B$ and $h=ab$ with $a\in A-C$, $b\in
B-C$. Furthermore $g^n\in C$ for some $n\in\mathbb N$ and $g^n$
corresponds to a power of the fiber.

\item $M_B$ has trivial JSJ-decomposition, $g\in C$ and $g$
corresponds neither to the fiber of $f'_Q$ of $Q^3$ nor to the fiber
of $M_B$ (if $M_B$ is Seifert).

\end{enumerate}\end{Lemma}

\noindent{\em Proof} We study the action of $G$ on the Bass-Serre
tree corresponding to the splitting $G=A*_CB$. We first investigate
the structure of the trees $T_g$ of elliptic elements of $G$. After
conjugation we can assume that $g\in A$ or $g\in B$.

If $g\in A-C$ then the tree consists of the vertex $x$ fixed under
the action of $A$
and all edges emanating from $x$. This is clear since in this case
$g^2$ corresponds to a power of the fibre of $Q^3$ and does therefore
not correspond to
the fibre of the piece $M_{B_1}$ of $M_B$ that contains $T$ (if
$M_{B_1}$ is Seifert).

If $g\in B$ then either
$T_g=y$ where $y$ is the vertex fixed under $B$ or after conjugation
$g^n\in C$ for some $n\in\mathbb N$. If $g\in C$, $M_{B_1}$ is
Seifert and $g$ is in $G$ conjugate to a power of $f_{B_1}$ we
conjugate $g$ such that $g=f_{B_1}^k$ for some $k\in\mathbb Z$. Note
that for $a\in A-C$ we get that $af_{B_1}a^{-1}\notin\langle
f_{B_1}\rangle$ since we assume that $f_{B_1}\neq f_A$. Note further
that $bcb^{-1}\notin C$ for all $b\in B-C$ and $c\in C-\langle
f_{B_1}\rangle$ since we assume that $M_B$ is not homeomorphic to
$Q^3$. This follows from the annulus theorem, the fact that all
properly immersed annuli can be homotoped into a Seifert piece and
since properly immersed annuli in Seifert pieces are homotopic to
vertical annuli. Therefore all essential properly immersed annuli in
$M_B$ with boundary in $T$ can be homotoped to be vertical in
$M_{B_1}$. Together this implies that $T_g$ must be contained in the
2-neighborhood of $y$.

As in the proof of Lemma~\ref{amalgam} we have to distinguish the
cases that $T_g\cap T_h\neq\emptyset$ and that $T_g\cap h T_g\neq
\emptyset$.

If $T_g\cap T_h\neq\emptyset$ we can assume that either $g$ or $h$ is
conjugate to an element of $A-C$, otherwise $\langle g,h\rangle$ lies
in the kernel of the quotient map $G\to G/N_G(B)\cong \mathbb Z_2$
which implies that $\langle g,h\rangle\neq G$. If $g$ and $h$ are
conjugate to elements of $A-C$ then we can argue as in the proof of
Lemma~\ref{amalgam} to obtain situation (iii) of 1. We are left with
the case
that $g\in A-C$ and that $h$ is conjugate to an element of $B$. The
structure of
the trees $T_g$ and $T_h$ guarantee that there exists a vertex $z$
with $hz=z$ such that either $d(x,z)=1$ or $d(x,z)=3$ where $x$ is
chosen such that $Ax=x$.

If $d(x,z)=1$ we have after conjugation with an element of $A$ that
$h\in B$, since further $g^2\in C$ we are in situation (i) of case
$1$.

If $d(x,z)=3$ then a power $h^k$ of
$h$ must be conjugate to an element of $C$ since it fixes an edge of
the Bass-Serre tree.
Note that this implies that $h$ is conjugate to
an element of $B_1$, the subgroup of $B$ corresponding to the piece
$M_{B_1}$ of the JSJ-decomposition of $M_B$ containing the torus $T$.
This is clear since otherwise this $h^k$ lies in two different edge
groups that are incident with $M_{B_1}$ which implies that $h^k$ is
conjugate to a power of the fibre of $M_{B_1}$, in particular
$M_{B_1}$ is Seifert. Such element however have no roots outside
${B_1}$.

We now study the action of $G$ on the Bass-Serre tree corresponding
to the splitting of $G$ that corresponds to the full
JSJ-decomposition of $M^3$. The trees $T_g$
and $T_h$ are now defined with respect to the new action of $G$.
Since $g\in A$ and a power of $h$ is conjugate to an element of $C$
it follows that $g$ and $h$ again act with fixed points. Since the
subgroup
generated by $g$ and $h$ is not free it follows that $T_g\cap T_h\neq
\emptyset$. As in the proof of Lemma~\ref{HNN1} we conclude that the
underlying graph
of the JSJ-decomposition must be a tree since the elliptic elements
$g$ and $h$ generate $G$. It is also clear that there is no other
piece in the JSJ-decomposition of $M^3$ that is homeomorphic $Q^3$,
otherwise we could define a surjective homomorphism $\phi:G\to
\mathbb Z_2$ whose kernel contains
$g$ and $h$.

It is easy to see that the trees $T_g$ and $T_h$ have the same
structure as before since we assume that $M_B$ is not homeomorphic to
$Q^3$. Thus $T_g$ lies in the $1$-neighborhood of $\bar x$ where
$\bar x$ is chosen such that $\bar x=A\bar x=g\bar x$ and $T_h$ lies
in the $2$-neighborhood of $\bar z$ where $\bar z$ is such that
$h\bar z=\bar z$ and $\bar z$ is fixed by a conjugate of $B_1$.

It follows that there exists a vertex $\bar y\in T_g\cap T_h$ with
$d(\bar x,\bar y)\le 1$ and $d(\bar z,\bar y)\le 2$. In particular we
have $d(\bar x,\bar z)\le 3$. Since $\bar x$ and $\bar z$ map onto
vertices of distance $1$ in the quotient tree (the tree underlying
the JSJ-decomposition and the corresponding graph of groups) we must
have that $d(\bar x,\bar z)=1$ or $d(\bar x,\bar z)=3$. In the first
case we argue as before. If
$d(\bar x,\bar z)=3$ then $d(\bar z,\bar y)=2$ and we denote the
midpoint of $[\bar y,\bar z]$ by $p$.
Since $h$ fixes $\bar z$ there must be a power of $h$ that fixes
$[\bar y,\bar z]$. We first
show that $[\bar y,\bar z]$ maps onto the edge of the JSJ-graph that
corresponds to $T$,
i.e. that $p$ is fixed by a conjugate of $A$. Suppose that $p$
corresponds to a piece $M'$ of the JSJ-decomposition of $M^3$ that is
different from $M_A$. It is clear that $M'$ is Seifert since $h^n$
lies in distinct conjugates of a peripheral subgroup of $\pi_1(M')$.
It further follows that $h^n$ corresponds to a power of $f_{M'}$ since
$M'$ is not homeomorphic to $Q^3$. This however yields a
contradiction since then $h$ cannot be conjugate to an element of
$C$. This implies that $\langle g,h\rangle \subset \langle
A,B_1\rangle$, in particular the
JSJ-decomposition of $M_B$ must consist of $M_{B_1}$ alone. Thus we
have $M_B=M_{B_1}$, $B_1=B$ and also $\bar x=x$, $\bar y=y$ and $\bar
z=z$.

We now look at the graph $[x,z]=[x,y]\cup [y,z]$. After conjugation
by an element of $A$ we can assume that $By=y$. Note that we can
consider the
segment $[y,z]=[y,p]\cup [p,z]$ as a single edge since $\Stab
[y,z]\cap\langle h\rangle =\Stab [y,p]\cap\langle h\rangle =\Stab
[p,z]\cap\langle h\rangle $. We
define $G_x=\langle g\rangle$, $G_y=\langle g^2,h^n\rangle$ where $n$
is chosen minimal such that $h^n$ fixes $[y,z]$ and $G_z=\langle
h\rangle$. We can assume that $n\ge 2$, otherwise $h$ fixes $y$ and
we are in situation (i) of 1. In particular we can assume that $h$ is
a root of the fibre of some conjugate of $B$ and that $h^n$ is the
fibre of $\Stab z\cong B$. This implies that $h^n$ (considered as an
element of $\Stab y$) has non-trivial and even intersection number
with the fibre of $B=\Stab
y$ since for any $c\in C$ and $a\in A-C$ the elements $c$ and
$aca^{-1}$ differ by an even multiple of the fibre $f_A$ (considered
as the Seifert space over $D(2,2)$). Since $f_A$ has by assumption
non-trivial intersection number (on $T$)
with $f_B$ we have that $af_Ba^{-1}$ has even intersection number
with $f_B$ for
all $a\in A-C$. We can exclude the situation that $M_B$ is a Seifert
space over $D(p,q)$ and that $g^2$ (a power of the fibre of $A$) has
intersection number $1$
with the fibre of $M_B$ since this puts us into case 2 of
Theorem~\ref{JSJ}. It then follows from Lemma~\ref{top10} that
$G_y=\langle
g^2,h^n\rangle$ is free in $g^2$ and $h^n$. This implies that $[x,y]$
and $[y,p]$
are $G_y$-inequivalent since otherwise $a^2$ and a conjugate of $h^n$
would commute which cannot happen if $\langle g^2,h^n\rangle$ is free
in $g^2$ and $h^n$. It is further clear that $G_y\cap \Stab
[x,y]=\langle g^2\rangle$ and that
$G_y\cap\Stab [y,p]=\langle h^n\rangle$ since primitive elements
generated maximal Abelian subgroups in free groups. By
Proposition~\ref{controlled} (2) we
therefore know that the induced splitting of $\langle g,h\rangle$ has
cyclic edge
groups which implies that $\langle g,h\rangle\neq G$.

\medskip\noindent If $T_g\cap hT_g\neq\emptyset$ we distinguish the
cases that $g$ lies (after conjugation) in $A-C$ or in $B$. If $g\in
A-C$ then the
structure of $T_g$ allows us argue as in the proof of
Lemma~\ref{amalgam} to get
in situation (ii) of case 1.

Suppose that $g\in B$. If $g\in C$ and $g$ has a root that is
conjugate to an element of $A-C$ we can clearly replace $g$ by this
element and argue as before. We can therefore assume that $g$ is not
conjugate to a power of $f_Q'$. If $T_g$ consists of a single vertex
then we argue as before that $\langle g,h\rangle$ fixes this vertex
which implies that $\langle g,h\rangle\neq G$. It follows that (after
conjugation) we can assume that $g^k\in C$ for some $k\in\mathbb Z$
and argue as before to see that
$g\in B_1$ where $B_1$ is as above. Again we assume that $g$ is a
root of the fibre of $B_1$ if $B_1$ is Seifert and $g$ is in $G$
conjugate to a root of the fibre of $B_1$. We can
assume that the JSJ-decomposition of $M^3$ has no piece that is
homeomorphic to $Q^3$ besides $M_A$, otherwise we could define a
surjective homomorphism from $G$ to $\mathbb Z_2* \mathbb Z_2$ such
that $g$ lies in the kernel which clearly implies that $\langle
g,h\rangle\neq G$. A similar argument shows that the JSJ-graph must
be a tree since we could otherwise define a surjective homomorphism
from $G$ onto $\mathbb Z_2* \mathbb Z$ with $g$ in the kernel.

We now study the action on the Bass-Serre tree associated to the
splitting that corresponds to the full JSJ-decomposition of $M^3$. It
is clear that $g$ again acts with a fixed point. By the remark after
Proposition~\ref{controlled} we can (after replacing $h$ by$g^kh$ for
some $k$) that either $T_g\cap T_h\neq \emptyset$ or $T_g\cap
hT_g\neq\emptyset$. If $T_g\cap T_h\neq\emptyset$ we are back in the
first case.

Suppose that $T_g\cap hT_g\neq\emptyset$. Let $y$ be the vertex fixed
under the action of $B_1$. As before we see that $T_g$ is contained
in the $2$-neighborhood of $y$. Since $M_B$ contains no piece that is
homeomorphic $Q^3$
it is further clear that for any vertex $v\in T_g$ with $d(y,v)=2$ we
see as before that $[y,v]$ maps onto the edge of the
JSJ-decomposition that corresponds
to $T$. We also have that $y$ is the only vertex of $T_g$ that can be
of valence
greater than $2$ since $C$ is of index $2$ in $A$. Choose two
vertices $p,q\in T_g$ such that $hp=q$. Since $p$ and $q$ project
onto the same vertex of
the JSJ-tree it follows that $d(p,q)$ is even. If $d(p,q)=0$ then $h$
acts with a fixed point and we are in the case before. We have to
deal with the cases that $d(p,q)=2$ and that $d(p,q)=4$.

If $d(p,q)=2$ then either $y\in\{p,q\}$ or $d(y,p)=d(y,q)=1$. In both
cases we have that $h$ lies in the subgroup of $G$ that is generated
by $B_1$ and the vertex group that corresponds to the vertex that is
the midpoint of $[p,q]$ in the first case or that corresponds to $p$
and $q$ in the second case. If $g$ and
$h$ generate $G$ this implies that this vertex corresponds to $A$ and
that the JSJ-decomposition of $M_B$ is trivial. In both cases we see
that after conjugation with an element of $B_1=B$ we have that the
edge fixed by $C$ lies in
$[p,q]$, it follows that $g^n\in C$ for some $n\in\mathbb N$. It is
further clear
that $h=ab$ for some $a\in A$ and $b\in B$ since $ap=q$ for some
$a\in A$. If $g\notin C$, i.e. $n\ge 2$ we have that $M_B$ is Seifert
fibered and that $g$ is
a root of the fibre. This puts us into situation 2 or 3.

Suppose now that $d(p,q)=4$. It is clear that $M_{B_1}$ is a Seifert
piece and that $g$ is a root of the fibre $f_{B_1}$ since otherwise
no power of $g$ could fix two edges emanating at $y$ (we assume that
$M_B$ is not homeomorphic to $Q^3$). As before we argue that
$af_Ba^{-1}$ has non-trivial even intersection number with $f_B$ for
all $a\in A-C$. Now $g^n$ and $hg^nh^{-1}$ fix $q$ and in $\Stab
y\cong B$ they correspond to boundary curves with even intersection
number
with the fibre that lie in different conjugacy classes of the
peripheral subgroup of $B$ corresponding to $T$. By Lemma~\ref{top10}
this implies that $\langle g^n,hg^nh^{-1}\rangle$ is free in $g^n$
and $hg^nh^{-1}$. We define $G_q=\langle g^n,hg^nh^{-1}\rangle$,
$G_p=\langle g^n,h^{-1}g^nh\rangle$ and $G_y=\langle g\rangle$. The
freeness of $G_q$ and $G_q$ and the minimality of $n$ guarantee as
before that $\Stab [y,q]\cap G_q=\Stab [y,q]\cap G_y$, that $\Stab
[y,p]\cap G_y=\Stab [y,p]\cap G_p$ and that $[y,q]$ and
$h[y,p]=[hy,q]$ are $G_q$-inequivalent. We can further
assume that $[y,p] $ and $[y,p]$ are $G_y$-inequivalent since
otherwise we had that $g^kh$ acts with a fixed point for some
$k\in\mathbb N$ which puts us into the first situation. Now by
Proposition~\ref{controlled} (3) we have that the induced splitting
of $\langle g,h\rangle$ has cyclic edge groups, i.e. $\langle
g,h\rangle\neq G$.~\hfill$\Box$

\medskip We conclude the proof of Theorem~\ref{JSJ}. We have the
different cases of
Lemma~\ref{specialcase}. In case (i)-(iii) of situation 1 the
argument is
the same as in the case without pieces of type $Q^3$. We therefore
only have to investigate
the cases that $h=ab$ and that either

\noindent (a) $g\in C$ and that $g$ corresponds to neither the fibre
$f'_Q$ of $Q^3$ nor
to the fibre of $M_B$ (if $M_B$ is Seifert) or

\noindent (b) $g\in B$ where $M_B$ is Seifert and $g$ is a root of
the fibre $f_B$.

\medskip
\noindent (a) We can assume that $g$ does not correspond to $f_Q$,
the fibre in the
fibration of $A$ as the orientable surface bundle over the M\"obius
band, since we can then argue as in the case without $Q^3$. We can
further assume
that $a^{-1}ga$ does not correspond to the fibre of $M_{B_1}$ where
$M_{B_1}$ is
the piece of $M_B$ containing $T$ since otherwise the proof of
Lemma~\ref{specialcase}
produces situation (b). The same arguments as in the case
without pieces of type $Q^3$ show that $\langle g,h\rangle\neq G$ if
$\langle g,h^{-1}gh\rangle\subset B$ is free. Also as in the case
without $Q^3$ we see that
$M_B$ must have a trivial JSJ-decomposition. Note, that $g\in C$ and
therefore
$a^{-1}ga=\bar g\in C$ since
$C$ is normal in $A$, therefore $h^{-1}gh=b^{-1}\bar gb\in B$.

\medskip Suppose that $M_B$ is a hyperbolic piece. Now $\langle
g,b^{-1}\bar gb\rangle$
cannot be Abelian since it does not lie in a parabolic subgroup.
Since we assume that $g$ is
neither a power of $f_Q$ not of $f_Q'$ it follows that $g\neq \bar
g^{\pm 1}$ which implies that $g$ and $b^{-1}\bar gb$ are not
conjugate in $B$ since $C$ is malnormal in $B$. By Lemma~\ref{hyp1}
we know that the covering of $M_B$ corresponding to $\langle
g,b^{-1}\bar gb\rangle$ must be homeomorphic to the exterior of a two
bridge link. This covering can however not be regular since $g$ and
$\bar g$ are not conjugate in $B$. This puts us into case~10 of
Theorem~\ref{JSJ}.

\medskip
Suppose that $M_B$ is Seifert. If $g$ and $\bar g$ have intersection
number
greater than $1$ with $f_B$ then $\langle g,b^{-1}\bar gb\rangle$ is
free by Lemma~\ref{top10}
and it follows that $\langle g,h\rangle\neq G$. The same holds if
$M_B$ is not a Seifert
piece over a base space of type $D(2,n)$. It remains to check the
case that the base space
of $M_B$ is of type $D(2,n)$ and that either $g$ or $\bar g$ has
intersection number one
with $f_B$. We can assume that the intersection number of $f'_Q$ with
$f_B$ is not $\pm 1$
since we are otherwise in case~2 of Theorem~\ref{JSJ}. W.l.o.g. we
can assume that $g$ has
intersection number $1$ with $f_B$. We use the notation as in the
beginning of
section~\ref{toplem}. We write $g=(xy)^n{f'_Q}^m$, it follows that
$\bar g^{-1}=(xy)^n{f'_Q}^{-m-2n}$, i.e. that $g$ and $\bar g$ differ
by $2(n+m)f'_Q$.
Since we assume that $f'_Q$ has at least intersection number $2$ with
$f_B$ this implies
that $\bar g$ has intersection number at least $4(n+m)-1$ with $f_B$.
It follows that the
intersection number of $\bar g$ with $f_B$ is greater than $3$ unless
$|n+m|\le 1$ in
which case we argue as before since $\langle g,b^{-1}\bar gb\rangle$
is
free by Lemma~\ref{top10}. If $n+m=0$ then $g$ corresponds to the
fibre $f_Q$ which
we have already excluded. It remains to check the case that
$|n+m|=1$. W.l.o.g. we can
assume that $n+m=1$, i.e. $m=-n+1$. The only case we have to deal
with is that
$g=(xy)^n{f'_Q}^{-n+1}$ has intersection number $1$ with $f_B$ and
that
$\bar g^{-1}=(xy)^n{f'_Q}^{-n-1}$ has (algebraic) intersection number
$-3$ with $f_B$.
This however implies that $(xy)^n{f'_Q}^{-n}$ has intersection number
$1$ ($-1$) with $f_B$.
This clearly implies that $n=\pm 1$, i.e. that $g$ corresponds to
$f_Q$ which we have
already excluded.

\medskip
\noindent (b) If the base space of $M_B$ is of type $D(p,q)$ or
$A(p)$ and the
intersection number of the fibres are one then we are in case~2 of
Theorem~\ref{JSJ},
otherwise the fundamental group cannot be generated by $g$ and $h$ by
Lemma~\ref{Q^3&S}.

If the base space is not of one of the above types then we consider
the quotient of $G$
by $N_G(C)$. The quotient map $\phi:G\to G/N_G(C)$ maps $G$ onto the
free product
$A/C*B/N_B(C)$. It is clear that $A/N_A(C)=A/C\cong\mathbb Z_2$ and
that $B/N_B(C)$ is
a orbifold group that is not generated by an elliptic element.
Furthermore $\phi(g)$ is
elliptic in the quotient orbifold group. The proof of Grusko's
theorem implies that
$\{\phi(g),\phi(h)\}$ is Nielsen-equivalent to a pair $\{\phi(g),\bar
h\}$ such that
$\bar h\in A/C\cup B/N_B(C)$. Such $\phi(g)$ and $\bar h$ however
cannot generate
$A/C*B/N_B(C)$ since $\phi(g)$ does not generate $B/N_B(C)$. It
follows that $\phi(g)$
and $\phi(h)$ cannot generate $A/C*B/N_B(C)$. Thus $g$ and $h$ do not
generate $G$.
\hfill$\Box$

\bigskip

\section{Heegaard genus}

In this section we deduce Theorem~\ref{heegrank} from
Theorem~\ref{JSJ}. We will
use the work of T. Kobayashi ([Ko1,2,3]).

In [Ko2], there is a complete list of closed orientable irreducible
3-manifolds with a genus 2 Heegaard splitting and with a non-trivial
JSJ-decomposition. Later
on, T. Kobayashi showed in [Ko3,\S3] how to extend the results of
[Ko1,2] in the
context of Heegaard splittings of irreducible 3-manifolds with
incompressible toral boundary. These results show that, except the
3-manifolds
given in Theorem~\ref{heegrank}, all the other ones
described in Theorem~\ref{JSJ} have Heegaard genus 2.

The fact that an irreducible Heegaard genus two splitting is strongly
irreducible,
implies the following
key lemma  (cf. [Mor, Lemma 1.1], [MSa, Lemma 2.2] or [RuS2,\S6]) :

\begin{Lemma} Let $M^3$ be a compact orientable irreducible
3-manifold with incompressible toral boundary and which admits a
Heegaard
decomposition
  $(V_1,V_2)$ of genus two. If $M$ has a non-trivial JSJ
decomposition, then
the JSJ family of tori  $\Sigma$ can be isotoped so that it intersects
$V_1$ and $V_2$ in a collection of essential annuli whose boundaries
are
also essential on
the tori $\Sigma$.\hfill$\Box$
\end{Lemma}

A proof of this lemma follows from the sweep-out of $M^3$ by the
Heegaard
surface $F = \partial V_1 = \partial V_2$ as described in [RuS1].

\

Now the description of Heegaard genus two, compact, orientable,
irreducible
3-manifolds with incompressible toral boundary follows, like in the
closed
case, from a
careful analysis of the possible intersections of the JSJ family
$\Sigma$
and the genus two
compression-bodies
$V_1$ and $V_2$ of the Heegaard decomposition (cf. [Ko2,\S3], [Mor,
Lemma 1.5],
[MSa, Theorem 2.1 and Lemma 2.2], [RuS2,\S5 and \S6]).

For example if $M^3$ has one boundary component, case 1) of
Theorem~\ref{JSJ} corresponds to
case 3-b) in the proof of Theorem 2.1 in [MSa]. In the same way,
case 2)
(respectively
cases 3) and 4) ) of Theorem~\ref{JSJ}
corresponds to case 2-b) (respectively 3-a) and 1) ) of the proof of
Theorem 2.1 in [MSa].
The analysis when $M^3$ has two boundary components is similar for
these
cases.

Case 6) of Theorem~\ref{JSJ} is explained in [Ko3, Lemma 6.1].

The cases 7), 8) and 9) of of Theorem~\ref{JSJ} corresponds to the
case
where the JSJ
family $\Sigma$  contains only non-separating tori. The analysis is
the same
as the one carried for the closed case in [Ko1, Thm 2] (cf. also
[RuS2,\S6,
cases 2) and 3)).

\
One can also deduce the
classification of
Heegaard genus 2 compact irreducible and $\partial$-irreducible
3-manifolds, with
non-empty toral boundary  and non-trivial JSJ-decomposition, from the
closed case by  using
the following lemma. This lemma is a direct consequence
of the work of Y. Rieck and E. Sedgwick [RiS1], [RiS2].

\begin{Lemma} Let $M^3$ be a compact orientable irreducible
3-manifold with incompressible toral boundary. Let $T^2\subset
\partial M^3$ be a boundary component. If infinitely many Dehn fillings on
$T^2$ yield a
manifold with a Heegaard splitting of genus 2, then $M^3$ has a
Heegaard splitting of genus 2. \end{Lemma}

\noindent {\em Proof} By the work of Y. Rieck and E. Sedgwick [RiS2, Cor.
6.6] there are
infinitely many Dehn fillings on $T^2$ such that the core $\gamma$ of
the attached solid torus can be isotoped into the genus 2 Heegaard
surface $\Sigma^2$ of the resulting manifold. Since a Heegaard genus two
splitting is
strongly irreducible, by [RiS1, Thm.4.6]  either:
\begin{enumerate} \item
 $\Sigma^2$ can be
isotoped to a Heegaard surface of $M^3$; or \item the surface
$\Sigma^2 \cap M^3$
is incompressible and $\partial$-incompressible in $M^3$ or compresses to
such an
essential surface.
\end{enumerate}

In the laste case, the filling slope is the boundary slope of an
essential surface, with all its boundary in $T^2$. By Hatcher's
finiteness result (cf. [Hat]), there are only finitely many such
slopes on $T^2$.

\hfill$\Box$

\section{Two-generated subgroups of 3-manifold fundamental groups}

In this section we prove Corollaries 4 to 7 which are direct
consequence of the proof
of Theorem~\ref{JSJ}.

\

\noindent The {\em proof of Corollary~\ref{twogen}} follows
immediately from the following Lemma:

\begin{Lemma}
Let $M^3$ be a compact orientable 3-manifold and let $U$ be a
2-generated
subgroup of $\pi_1(M^3)$ that is neither free or free abelian. If $U$
is
of infinite index in $\pi_1(M^3)$, then  $U$ is either the
fundamental group of
a Seifert fibered manifold, or of a complete hyperbolic manifold with
finite volume,
  or of one of the manifolds described in Theorem~\ref{JSJ}.
\end{Lemma}

\noindent {\em Proof}
Let $p: \hat M^3 \to M^3$ be the covering
of $M^3$ with fundamental group $\pi_1(\hat M^3) = U$. By Scott's
compact core theorem,
there is a compact submanifold $\hat N^3 \subset \hat M^3$ such that
the
inclusion
map $\hat N^3 \hookrightarrow \hat M^3$ induces an isomorphism between
$\pi_1(\hat N^3)$ and $U$.

If $U = U_1 \star U_2$ is a non-trivial free product, then
one  factor, say $U_1$, is a finite cyclic group since $U$ is not a
free group.
In this case, since $M^3$
is orientable, by Epstein's theorem [Ep] (cf.
[He, Thm.9.8]) $M^3 = M^{3}_{1} \# R^3$ where $R^3$ is closed and
orientable,
$\pi_1(R^3)$ is finite. Moreover, after conjugation, one may assume
that $U_1$ is  a
subgroup of $\pi_1(R^3)$.
Since the covering $p: \hat M^3 \to M^3$ is of infinite index,
$p^{-1}(R^3)$ would lift to
infinitely many connected sum factors of $\hat M^3$ with non-trivial
fundamental group $U_1$,
contradicting the fact that $\pi_1(\hat M^3) = U$ is of rank two.

Therefore, we can always assume that $U$ is freely indecomposable.
Since $U = \pi_1(\hat N^3)$ is not infinite cyclic,
according to [He, Lemma 10.1] there is a compact orientable and
irreducible 3-manifold
$N_{0}^3$ with fundamental group $\pi_1(N_{0}^3) = U$. It is the
Poincar\'e
associate of $\hat N^3$, i.e. the only non-simply connected prime
factor of the 3-manifold
obtained by capping off the 2-spheres in $\partial \hat N^3$
by 3-balls.

Since $U$ is of infinite index in $\pi_1(M^3)$, its cohomological
dimension is smaller or equal to two. Hence it cannot be the
fundamental group of a
closed irreducible 3-manifold. Therefore the compact, orientable,
irreducible
3-manifold
$N_{0}^3$ has a non-empty incompressible boundary, since $U$ is
not a free product. Hence $N_{0}^3$
is a Haken 3-manifold with incompressible boundary. According to
[JS2, Lemma
5.4], all
components of $\partial N_{0}^3$ are tori.

If $N_{0}^3$ is atoroidal (i.e. has a trivial JSJ decomposition), it
follows
from Thurston's hyperbolization theorem that it is either a Seifert
fibered 3-manifold
or a complete hyperbolic 3-manifold with finite volume.

If $N_{0}^3$ has a non-trivial JSJ-decomposition then $N_{0}^3$ is
one of the manifold described in Theorem~\ref{JSJ}.
\hfill$\Box$

\

\noindent {\em Proof of Corollary~\ref{twoperiph}}
The following can be extracted from the proof of Theorem~\ref{JSJ}
but a direct argument is just as easy. We first show that $M^3$ must
have a trivial
JSJ-decomposition if $\pi_1(M^3)$ is generated by two peripheral
elements. Suppose that $\pi_1(M^3)=\langle g,h\rangle$ where $g$ and
$h$ are peripheral. It is clear that $g$ and $h$ act as elliptic
elements on the Bass-Serre tree asscociated to the splitting of
$\pi_1(M^3)$ corresponding to the JSJ-decomposition. It then follows
from \cite{KW} that powers of $g$ and $h$ must have a common fixed
point. Since $g$ and $h$ are peripheral they are not proper roots of
the fiber of some Seifert piece in the JSJ-decomposition. Thus $g$
($h$) fixes every point in the Bass-Serre tree that is fixed by a
power of $g$ ($h$). It follows that $g$ and $h$ have a common fixed
point, i.e. $\pi_1(M^3)=\langle g,h\rangle$ is conjugate to a vertex
stabilizer. It follows that the splitting and therefore the
JSJ-splitting of $M^3$ is trivial.

Since it has
an incompressible
boundary, according to [JS2, Lemma 5.4] all
components of $\partial M^3$ are tori. Thus by Thurston's
hyperbolization theorem
$M^3$ is either a Seifert fibered 3-manifold
or a complete hyperbolic 3-manifold with finite volume.
In the Seifert fibered case, the proof follows from Lemmas~\ref{top3}
and \ref{top6}.
In the hyperbolic case it follows from Lemma~\ref{hyp1}.
\hfill$\Box$

\

\noindent {\em Proof of Corollary~\ref{oneperiph}}
Let $M^3$ be a compact orientable, irreducible
3-manifold with incompressible boundary. If $M^3$ is not hyperbolic
and
$\pi_1(M^3)$ is generated
by two elements one of which is peripheral then either $M^3$ has
Heegaard genus $2$ or $M^3$ belongs to case 2) of
Theorem~\ref{heegrank}.

So we can assume that $M^3 = S\cup_{T} H$. Here $H$ is a hyperbolic
3-manifold with
$\pi_1(H)$ generated by a pair of elements with a  parabolic generator
belonging to
the parabolic subgroup associated to the boundary component $T$. The
Seifert manifold
$S$ has basis $D(p,q)$ or $A(p)$. The gluing map identifies the fibre
of
$S$ on $T$
with the simple closed curve corresponding to the parabolic generator
of
$\pi_1(H)$.

We will show that $H$ is homeomorphic to the exterior of a two bridge
link and that the parabolic generator corresponds to a meridian. It
then follows from \cite{Ko2} that $M^3$ has Heegaard genus two.

Suppose now that $\pi_1(M^3)=\langle g,h\rangle$ where $g$ is
peripheral. As before we see that $g$ acts with fixed point on the
associated Bass-Serre tree and it follows from \cite{KW} that we can
choose $h$ such that either $T_g\cap T_h\neq\emptyset$ or $T_g\cap
hT_g\neq\emptyset$. In particular the conclusion of
Lemma~\ref{amalgam} and Lemma~\ref{specialcase} hold with this
particular $g$ except possibly in the case that $h$ also acts with a
fixed point; here we might have to exchange $g$ and $h$. It now
suffices to observe that in the proof of Theorem~\ref{JSJ} either $H$
was the exterior of a 2-bridge link or $g$ was a proper root of the
fibre. The last case however is impossible if $g$ and $h$ were not
interchanged since a peripheral element cannot be a proper root of
the fibre. If $g$ and $h$ were exchanged then the proof of Theorem
\ref{JSJ} implies that $g$ is a root of the fibre and $\pi_1(H)$ is
conjugate to $\langle h,g^n\rangle$ where $g^n$ is a power of $g$
that lies in the conjugate of some edge group. In particular $g^n$ is
peripheral in $\pi_1(H)$. Since $h$ is also peripheral it follows
from Corollary~\ref{twoperiph} that $H$ is the exterior of a 2-bridge
link and that $h$ corresponds to a meridian.\hfill$\Box$

\

\noindent {\em Proof of Corollary~\ref{sat}}
Let $k \subset \mathbb S^3$ be a 2-generator satellite knot. It
follows from
Theorem~\ref{heegrank} that $k$ is tunnel number one or its exterior
$E(k)$
belongs to case 2) of Theorem~\ref{heegrank}.

It means that $E(k) = M_1\cup_{T} M_2$, where $M_1$ is a Seifert
manifold
and $M_2$ is
a hyperbolic 3-manifold.

The Seifert manifold $M_1$ has basis $D(p,q)$ or $A(p)$, hence it is
the
exterior of a $(p,q)$-torus knot in
$\mathbb S^3$,
or $k$ is a cable knot [Ja, Lemma IX.22]. But, by using the cyclic
surgery Theorem [CGLS],
Bleiler [Ble2]
has shown that  a 2-generator cable knot is an iterated
torus knot, contradicting the fact that $M_2$ is hyperbolic.

Therefore $M_1$ is the exterior of a torus knot and $\partial M_2 = T
\cup
\partial E(k)$
has two components. The splice decomposition  of the satellite knot $k$
shows that $M_2$ is homeomorphic to the exterior of a link
$L' = k'_0 \cup k'_1$ obtained in the  way described in [EN, Prop. 2.1].
The torus $T$
bounds in $\mathbb S^3$ a solid torus $V_0$ containing $k$ and $M_2$ is
homeomorphic
to the exterior of $k$ in $V_0$. By gluing a solid torus $V_1$ to $T=
\partial V_0$ in
such way that the prefered longuitude of the torus knot exterior on $T$ is
identified
with the boundary of a meridian of $V_1$ one gets $\mathbb S^3$.
Then the image of $k$  in $\mathbb S^3$  by this
desplicing operation is $k'_0$ and $L'$ is the link
formed by $k'_0$ and the core  $k'_1$ of the unknotted solid torus $V_1$.
It follows that $k'_1$ is unknotted in $\mathbb S^3$

The  hyperbolic 3-manifold $M_2$ has a fundamental group
$\pi_1(M_2)$ generated by a pair of elements with a single parabolic
generator $a$ belonging to
the parabolic subgroup associated to the boundary component $T$. The
gluing
map identifies
the fibre of $S$ on $T$ with the simple closed curve $\alpha \subset
T$
corresponding to the parabolic
generator $a$.

Because of the gluing instruction, this curve $\alpha$ intersects on
$T$
the meridian of the torus knot
and hence the meridian of $k'_1$
only once. Therefore a Dehn surgery along $k'_1$ with slope $\alpha$
yields
back $\mathbb S^3$.
Let $k_1$ be the core of this Dehn surgery and $k_0 \subset \mathbb
S^3$ be
the image of $k'_0$
after this Dehn surgery. Then the exterior of $L = k_0 \cup k_1$ is
homeomorphic to the exterior
of $L'$, hence to $M_2$.

Moreover, $\pi_1(M_2)$ is generated by two elements where one is a
meridian
of the component
$k_1$ of $L$. It follows that the fundamental group of the exterior
of the
other component
$k_0$ is cyclic, hence $k_0$ is also unknotted by the loop theorem.

Then it follows from [Ku] that $L$ has tunnel number one if and only
if it is a
2-bridge link. But in this
case $K$ has tunnel number one (cf. [MSa]).
\hfill$\Box$

\section{Involutions on 3-manifolds with a 2-generated fundamental
group}
  In this section we prove Corollaries \ref{branchcov} and \ref{geom}.

\subsection{2-fold branched coverings of homotopy spheres}

For the proof of Corollary \ref{branchcov} we distinguish the two
cases
that $M^3$ is geometric and that $M^3$ has a non-trivial
JSJ-decomposition.

\noindent {\bf 1. case: } $M$ is geometric.

Then it is either a Seifert fibred, a Sol or a hyperbolic
3-manifold with \allowbreak $\rank \pi_1(M^3)=2$.

In the Seifert fibred case, the proof  follows from the determination
of
rank 2 closed Seifert
3-manifolds [BZg] and the construction of the fibre preserving
Montesinos involution
on these manifolds which shows that they are 2-fold branched coverings of
$\mathbb S^3$
([Mon1],[Mon2, \S4.7]).

For Sol 3-manifolds, the proof follows from the determination
of rank 2 orientable torus bundles [TO,Lemma 1] (cf. also [Sa2]).

In the hyperbolic case, the proof follows from
the following proposition  which includes the
case with boundary. We will need it  when the JSJ-decomposition is not
trivial.
We recall that a n-times punctured homotopy
3-sphere is a 3-dimensional homotopy sphere minus the interior of $n$
disjoint
embedded 3-balls.

\begin{Proposition}\label{fakesphere} Let $M^3$ be a compact
orientable complete
hyperbolic 3-manifold with finite volume. If \allowbreak $\rank
\pi_1(M^3)=2$ and $\pi_1(M^3)$ is not abelian,
then $M^3$ is a 2-fold branched covering of a n-times punctured
homotopy
sphere with $n \leq 2$. \end{Proposition}

\noindent {\em Proof} Let $\alpha$ and $\beta$ two elements in
$\pi_1(M)$ that generate the group. The proof follows from the
following
lemma due to Jorgensen (cf.
[Th, chap.5, Prop.5.4.1 and 5.4.2] ):

\begin{Lemma}\label{invo} Any complete hyperbolic 3-manifolds $M^3$
whose
fundamental group is not abelian and generated by two elements
$\alpha$ and $\beta$
admits a non-free, orientation preserving isometry $\tau$  of order
2  which takes
$\alpha$ to $\alpha^{-1}$ and $\beta$ to $\beta^{-1}$. Moreover if
the two generators are
conjugate in $\pi_1(M^3)$, there is an orientation preserving,
isometric
$\mathbb Z^2 \oplus \mathbb Z^2$ action on $M^3$ generated by $\tau$,
together with
an involution $\rho$ which interchanges $\alpha$ and $\beta$.
\end{Lemma}

\noindent {\em Proof of Lemma \ref{invo}}
The involution $\tau$ is induced by conjugation by the rotation $t$
of angle $\pi$ about the the following line $\Delta_0$ in the
hyperbolic space
$\mathbb{H}^3$:

(i) $\Delta_0$ is the unique common perpendicular to the axis of
$\alpha$ and
$\beta$ if they are both loxodromic elements in $PSL_2(\mathbb C)$;

(ii) $\Delta_0$ is the unique line through the fixed points at
infinity of $\alpha$
and $\beta$, if they are both parabolic elements in $PSL_2(\mathbb
C)$,(the two fixed
points are distinct, because $\pi_1(M)$ is not abelian);

(iii) $\Delta_0$ is the unique perpendicular to the axis of $\alpha$
through the
fixed point at
infinity of $\beta$, if $\alpha$ is loxodromic and $\beta$ is
parabolic.

In particular, $t \alpha t^{-1} = \alpha^{-1}$ and $t \beta t^{-1} =
\beta^{-1}$. Therefore, if $w(\alpha,\beta) = 1$ is any relation in
the group $\pi_1(M)$, then $w(\alpha^{-1},\beta^{-1}) = 1$. It
follows that the rotation
$t$ conjugates the group $\pi_1(M)$ to itself, and hence induces a
non-free,
orientation preserving, isometric involution
$\tau$ on $M^3$ which takes
$\alpha$ to $\alpha^{-1}$ and $\beta$ to $\beta^{-1}$.

If $\alpha$ and $\beta$ are conjugated in $\pi_1(M^3)$, then they are
either
both loxodromic elements or both parabolic elements (with distinct
fixed point)
in $PSL_2(\mathbb C)$.
Then the two involutions $\rho$ and $\tau \circ \rho$ are induced by
conjugation by  rotations
$r_1$ and $r_2$ of angle $\pi$ about two  lines $\Delta_1$ and
$\Delta_2$
that are perpendicular to each other and perpendicular to the line
$\Delta_0$
in the hyperbolic space
$\mathbb{H}^3$:

(iv) If $\alpha$ and
$\beta$  are both loxodromic elements in $PSL_2(\mathbb C)$,
$\Delta_1$ and $\Delta_2$
intersects $\Delta_0$ at
the mid point between the axis $\Delta_{\alpha}$ and $\Delta_{\beta}$
of
$\alpha$ and $\beta$. Moreover  the planes  $(\Delta_0, \Delta_1)$
and
$(\Delta_0, \Delta_2)$ are the two bisector planes between the
two planes $(\Delta_0, \Delta_{\alpha})$ and $(\Delta_0,
\Delta_{\beta})$.

(v) If $\alpha$
and $\beta$ are both parabolic elements in $PSL_2(\mathbb C)$,
$\Delta_1$
and $\Delta_2$ intersect
$\Delta_0$ at the unique point $p \subset \Delta_0$ such that $d(p,
\alpha(p) = d(p, \beta(p)$.

In particular the two rotations $r_1$ and $r_2$ around $\Delta_1$ and
$\Delta_2$ commute
with the rotation $t$ around $\Delta_0$.

One of these rotations, say $r_1$  conjugates
$\alpha$ to $\beta$, while the second one $r_2 = t\circ r_1$
conjugates $\alpha$ to $\beta^{-1}$.
Therefore, if $w(\alpha,\beta) = 1$ is any relation in
the group $\pi_1(M)$, then $w(\beta, \alpha) = 1$. So the rotation
$r_1$ conjugates the group $\pi_1(M)$ to itself, and hence induces a
non-free,
orientation preserving, isometric involution
$\rho$ on $M^3$ which takes
$\alpha$ to $\beta^{-1}$ and commutes with the involution $\tau$.
\hfill$\Box$

\

We finish now the proof of Proposition \ref{fakesphere}.
Let $\Gamma$ be the subgroup of $PSL_2(\mathbb C)$ generated by
$\alpha$, $\beta$
and $t$. Then $\pi_1(M)$ is a subgroup of index at most 2 in
$\Gamma$.
Hence $\Gamma$
is a discrete cocompact subgroup of $PSL_2(\mathbb C)$ and ${\cal
O}=\mathbb H^3/\Gamma$ is the compact orientable hyperbolic 3-orbifold
$\cal O = M/\tau$ obtained as the quotient of $M^3$ by the
orientation
preserving isometric involution $\tau$. The orbifold fundamental
group
of $\cal O$ is $\Gamma$.

Since $\Gamma$ is generated by the three orientation preserving
isometric involutions $t \alpha$, $t \beta$ and $t$ of the
hyperbolic space $\mathbb H^3$, the fundamental group of the
underlying space $\vert \cal O \vert$
of $\cal O$ is trivial (cf. [Th, Chap.13]). In particular $\vert \cal
O \vert$
is a n-times punctured homotopy sphere , and $M^3$ is a 2-fold
branched covering
of it.

Since $\rank \pi_1(M^3)=2$, $M^3$ has at most two boundary components
and they are tori.
Therefore $\partial \cal O$ has at most two boundary components and
they are {\em pillows}
(i.e. 2-spheres with 4 branching points of order 2). Thus
$\vert \cal O \vert$ is a n-times punctured homotopy sphere with $n
\leq 2$.
In particular the restriction of the involution $\tau$ on each torus
boundary component is
conjugated by an isotopy to the Weierstrass involution.
\hfill$\Box$

\begin{Remark}

Except in the case where $M^3$ has Heegaard genus 2, we cannot
prove that $\vert \cal O \vert$ is the  true sphere $\mathbb S^3$,
the true ball $\mathbb B^3$ or the product $S^2 \times [0,1]$.
\end{Remark}

\noindent {\bf 2. case: } $M^3$ has a non-trivial JSJ-decomposition.

Then either $M^3$has Heegaard genus two or belongs to one of the
cases 2) to 4)
described in Theorem~\ref{heegrank}.

\

If $M^3$ has Heegaard genus two, then it is a 2-fold branched
covering of $\mathbb S^3$
by Birman and Hilden [BH].

\

If $M^3$ belongs to case 2), $M^3 = S \cup_{T} H$, where $H$ is a
hyperbolic 3-manifold
with a 2-generated fundamental group and $S$ is a Seifert 3-manifold
over $D(p,q)$.

By Proposition~\ref{fakesphere} $H^3$ is a 2-fold branched covering
of a homotopy ball
and the restriction of the covering involution to $T = \partial H$ is
conjugated by an isotopy
to the Weierstrass involution.

A Seifert manifold $S$ with basis $D(p,q)$ admits an orientation
preserving Montesinos involution with quotient a $3$-ball, which is
fiber-preserving and reverses the orientation of the fiber. The
restriction of this
involution to $T = \partial S$ is also conjugated by an isotopy to
the Weierstrass involution.

Since any gluing homeomorphism $f : T \to T$ is isotopic to one
which commutes with the Weierstrass involutions on $T$, one can glue
the
two involutions to get an involution on $M^3$ with quotient a
homotopy sphere.

\

If $M^3$ belongs to case 3), $M^3 = S_{1} \cup_{T} S_{2}$, where
$S_1$ is a Seifert 3-manifold
over $M\ddot o$ or $M\ddot o(p)$ and $S$ is a Seifert 3-manifold over
$D(2,2l+1)$.
Both sides admit
an orientation preserving Montesinos involution whose restriction to
the torus boundary is
conjugated by an isotopy to the Weierstrass involution. The same
argument as in case 2) shows
that $M^3$ is a 2-fold branched covering of $\mathbb S^3$.

\

If $M^3$ belongs to case 4), $M^3 = Q \cup_{T} H$, where $H$ is a
hyperbolic 3-manifold
and $Q$ is the orientable circle bundle over $M\ddot o$. The
difference here is that
we do not know
that the fundamental group of the hyperbolic piece is generated by
two elements. However,
the restrictions on $\pi_1(H)$ given at the end of $\S 3.3$ allow to
show:

\begin{Lemma}\label{irreg} In case 4) of Theorem~\ref{heegrank}, the
fundamental
group of the hyperbolic piece $H$ is generated by two conjugated
peripheral subgroups
$P_1 = \pi_1(\partial H)$ and $hP_{1}h^{-1}$, with $h\in
\pi_1(H)$.\hfill$\Box$
\end{Lemma}

By an hyperbolic Dehn filling argument we deduce from
Lemma~\ref{irreg} the following
corollary:

\begin{Corollary}\label{fakeball} In case 4) of
Theorem~\ref{heegrank}, the hyperbolic
piece $H$ is a 2-fold branched covering of a homotopy ball.
\end{Corollary}

\noindent {\em Proof of Corollary~\ref{fakeball}}
Let $\alpha \subset \partial H$ be a simple closed curve and let
$H(\alpha)$ be
the closed 3-manifold
obtained by gluing a solid torus along $\partial H$ in such way that
$\alpha$
is identified with the boundary of a meridian disk. It follows from
Lemma~\ref{irreg}
$\pi_1(H(\alpha))$ is generated by two conjugated elements.

By Thurston's hyperbolic Dehn filling theorem [Th] (cf.[BP],[CHK]), for
almost all simple
curves $\alpha \subset \partial H$, $H(\alpha)$ is  hyperbolic and
the core $k_{\alpha}$
of the Dehn filling is the shortest geodesic in $H(\alpha)$.

By Lemma~\ref{invo} there is an orientation preserving, isometric
$\mathbb Z^2 \oplus \mathbb Z^2$ action on $H(\alpha)$ generated by
an involution $\tau$
  which take each generator of $\pi_1(H(\alpha))$ to its inverse,
together with
an involution $\rho$ which interchanges the two generators. Moreover,
this
$\mathbb Z^2 \oplus \mathbb Z^2$ action preserves the shortest
geodesic $k_{\alpha}$.
Hence, it induces a $\mathbb Z^2 \oplus \mathbb Z^2$ action on the
exterior $H$ of
$k_{\alpha}$. Moreover, the involutions $\tau$  reverses the
orientation of
$k_{\alpha}$ since it takes each generator to its inverse. In
particular, the quotient of
an invariant tubular neighborhood of $k_{\alpha}$ by $\tau$ is a
3-ball.

By Proposition~\ref{fakesphere} the quotient of $H(\alpha)$ by the
involution $\tau$ is
a homotopy sphere, hence the quotient of $H$ by the restriction of
$\tau$ is a
homotopy ball.
\hfill$\Box$

\

Using Corollary~\ref{irreg} and the fiber preserving Montesinos
involution on $Q$, viewed
as a Seifert manifold over $D(2,2)$, one can conclude in case 4) like
in case 2).
This finishes the proof of Corollary~\ref{branchcov}.
\hfill$\Box$

\

Examples of closed hyperbolic 3-manifolds with a 2-generated
fundamental group
can be obtained by hyperbolic Dehn filling along a 2-generator knot
in $S^3$,
for examples a 2-bridge knot.

Using hyperbolic Dehn fillings of once punctured torus bundles with
pseudo-Anosov monodromy,
  A. Reid [Rei] has
been able to construct infinitely many closed 2- or 3-generator
hyperbolic
3-manifolds
which have a proper
finite sheeted cover with a 2-generated fundamental group. From his
construction,
he gave also examples of
closed Haken hyperbolic 3-manifolds which are neither a surface
bundle nor
double covered by a
surface bundle, but which have a finite cover which has 2-generated
fundamental group.

However, a 2-generator once punctured torus bundle has Heegaard genus
2 by
[TO,Lemma 1] (cf. also [Sa2]).

We give now an example of  hyperbolic genus 2 surface bundle
with a 2-generated fundamental group, due to Nielsen.

Let $F$ be a compact orientable surface with negative Euler
characteristic. We
say that a diffeomorphism $\phi \in Diff^{+}(F)$ fills up $\pi_1(F)$
if there is
an element $\gamma \in \pi_1(F)$ such that its orbit
$\{\phi_{\star}^{n}(\gamma)\}_{n\in \mathbb Z}$ generates $\pi_1(F)$.

It is an easy observation that the mapping torus of a diffeomorphism
$\phi \in Diff^{+}(F)$ that fills up $\pi_1(F)$ has a 2-generated
fundamental group.

The following example (Nielsen's example $\sharp 13$ in [Nie]; cf.
[Gil]) is a
pseudo-Anosov diffeomorphism of a one-punctured surface $F$ of genus 2
which fills up
$\pi_1(F)$. Hence its mapping torus $M_{\phi}$ is a complete
hyperbolic
3-manifold with finite
volume and one cusp.

The closed genus 2 surface bundle, obtained by Dehn filing $\partial
M_{\phi}$ along
the the closed curve $\partial F$, has still a pseudo-Anosov
monodromy and
hence has
a 2-generated fundamental group.
\

A diffeomorphism
$\phi \in Diff^{+}(F)$ is determined, up to isotopy, by its induced
action on
$\pi_1(F) = \langle a,b,c,d \rangle$, where the conjugacy class of
$[a,b][d,c]$ represents the simple loop $\partial F$.

Let  $\phi_{\star}: \pi_1(F) \to \pi_1(F)$ be the automorphism given
by:

$$
\phi_{\star}(a) =  c^{-1}a^{-1} \,\,;\,\, \phi_{\star}(b) =
b^{-1}a^{-1}
$$
$$
\phi_{\star}(c) = b^{-1}a^{-1}d  \, \, ; \,\, \phi_{\star}(d) =
c^{-1}.
$$

One easily verifies that :

$$ \phi_{\star}([a,b][c,d]) = (abc)^{-1}[a,b][c,d](abc)$$.

That is the necessary and sufficient condition for $\phi_{\star}$ to
be
induced by a
diffeomorphism $\phi \in Diff^{+}(F)$.

Moreover, it is easy to check that $\pi_1(F)$ is generated by
$\langle b,
\phi_{\star}(b),
  \phi_{\star}^{2}(b), \phi_{\star}^{3}(b) \rangle$. Hence $\phi$
fills up
$\pi_1(F)$.

The Nielsen 's classification of diffeomorphisms of surfaces [Nie],
based
on the Nielsen types of the lifts of the diffeomorphism to the unit
disk
allows to show that $\phi$ is pseudo-Anosov (cf. [Gil]).

\subsection{2-generator knots in 3-manifolds}

The proof of corollary \ref{geom} follows from Thurston's orbifold
theorem (cf.[BP],[CHK]) and
the following lemma:

\begin{Lemma}\label{invertibe} Let $k\subset M^3$ be a
2-generator knot in a closed, orientable,
irreducible 3-manifold $M^3$. If $k$ is not contained in a ball, then
$k$ is
strongly invertible: i.e.
there is an orientation preserving involution on $M^3$
that preserves $k$ while reversing its orientation.
\end{Lemma}

\noindent {\em Proof} Since $M^3$ is irreducible and $k\subset M$ is
not contained in a
ball, the exterior $E(k)$ of $k$ is irreducible. It follows from
Theorem~\ref{heegrank} that $E(k)$ has Heegaard genus two or
  belongs to the case 2) described in Theorem~\ref{heegrank}.

If $E(k)$ has Heegaard genus 2, then
the hyperelliptic involution on the genus 2 Heegaard surface extends
to
both sides
since the attaching curve for the 2-handle can be choosen, up to
isotopy,
to be invariant
under the hyperelliptic involution. This gives a non-free, orientation
preserving involution
$\tau$ on $E(k)$ whose restriction to $\partial E(k)$
is conjugated to the  Weierstrass involution with quotient a 2-sphere
with 4 branching points. In particular it
extends to an orientation preserving, non-free involution  $\bar \tau$
on $M$, preserving $k$ and with two fixed points
on $k$.

\

If $E(k)$ belongs to the case 2) described in Theorem~\ref{heegrank},
$E(k) = S \cup_{T} H$, where $H$ is a hyperbolic 3-manifold
with a 2-generated fundamental group and $S$ is a Seifert 3-manifold
over $D(p,q)$ or A(p).

By Proposition~\ref{fakesphere} $H$ is a 2-fold branched covering
of a homotopy ball or twice punctured homotopy sphere, according
whether
$\partial H$ has one or two components. Moreover
the restriction of the covering involution to each component $T
\subset
\partial H$ is
conjugated by an isotopy
to the Weierstrass involution.

A Seifert manifold $S$ with basis $D(p,q)$ (respectively  $A(p)$)
admits an
orientation
preserving Montesinos involution with quotient a $3$-ball
(respectively a product $S^2 \times [0,1]$),
which is
fiber-preserving and reverses the orientation of the fiber. The
restriction of this
involution to each component $T \subset \partial S$ is also
conjugated by
an isotopy to
the Weierstrass involution.

Then a gluing argument, analogous to the one used in the proof of
Corollary~\ref{branchcov}, shows that $E(k)$ admits an orientation
preserving,
non free involution $\tau$ whose whose restriction to $\partial E(k)$
is conjugated to the  Weierstrass involution. As above it
extends to an orientation preserving, non-free involution  $\bar \tau$
on $M$, preserving $k$ and with two fixed points
on $k$.
\hfill$\Box$

\

\noindent {\em Proof of corollary \ref{geom}} Since $M$ is
irreducible, by
the proof of the
Smith conjecture ([MB],[Wa3])
the orbifold $M/\bar{\tau}$ is irreducible. Then, by Thurston's
orbifold
theorem ([BP], CHK) either $M/\bar{\tau}$ is geometric, or $M$ has a
non-trivial $JSJ$-decomposition.~\hfill$\Box$

\section{2-generator knots in $\mathbb S^3$}

In this section we prove corollary 3, which extends to 2-generator
knots a result shown for tunnel number one knots by M.Scharlemann
([Scha]).

We recall that a knot
$k \subset \mathbb S^3$ is prime if there is no sphere $S^2$ that
meets $k$
transversally in two points and intersects the exterior $E(k) =
\mathbb{S}^{3} -
int(N(k))$ in an essential annulus.

A Conway sphere for a knot $k$ is sphere $S^2$ that meets $k$
transversally in
four points and gives in the exterior $E(k)$ of $k$ an essential
planar surface.
The knot $k$ is Conway irreducible (or doubly prime) if there is no
Conway
sphere for $k$.

The fact that 2-generator knots are prime follows already from [N]
(see also [W1]).
So in the remaining we assume hat the knot $k$ is prime. Then by [BS]
for a
prime knot $k \subset \mathbb S^3$ there is a finite characteristic
collection of
disjoint, non-parallel tori and Conway spheres such that:\ (i) the
collection of
tori corresponds to the JSJ-family of tori of the exterior $E(k)$ of
$k$; (ii) the
collection of tori and Conway spheres lifts to the JSJ-family of tori
of the
2-fold branched covering of~$k$.

In particular a knot $k\subset \mathbb S^3$ is prime and Conway
irreducible if and only if its
2-fold branched covering is irreducible and topologically atoroidal
(cf. [Ble1]).

Starting with a 2-generator knot $k$, we consider its
Bonahon-Siebenmann characteristic
collection $\cal C$ of tori and Conway spheres, and we assume that
$\cal C$ contains at least one Conway sphere. This Conway sphere will
avoid the JSJ
family of tori $\cal T$ in $\cal C$.

Let $X^{3} \subset E(k)$ be the closure of the connected component of
$E(k) - \cal T$
that contains $\partial E(k)$. Then the track of the Conway sphere in
$E(k)$ is an essential,
properly embedded four punctured sphere
$(\Sigma^{2}, \partial \Sigma^{2}) \hookrightarrow (X^{3}, \partial
E(k))$.

If $\rank(\pi_1(E(k)))=2$, it follows from Theorem~\ref{JSJ} that
either
$X^{3}$ is
Seifert fibered, or it is hyperbolic with $\rank(\pi_1(X^{3}))=2$.

If $X^{3}$ is Seifert fibered, the only essential surfaces in
$X^{3}$, with
non-empty boundary, are either saturated annuli or horizontal
surfaces
(i.e branched covering of the basis of the Seifert fibration).
Hence it is an annuIus
or it meets each boundary component of $\partial X^3$. In
both cases,
it cannot be a four punctured sphere with only meridional boundary
components (cf. also [GL, Lemma 3.1], since $X^3$ must be a cable space).

If $X^{3}$ is hyperbolic and $\rank(\pi_1(X^{3}))=2$, one uses the
following well-known lemma:

\begin{Lemma} Let $M^{3}$ be a compact, orientable, irreducible and
atoroidal 3-manifold
with empty or incompressible boundary. If $\rank(\pi_1(M^{3}))=2$,
then $M^{3}$ does not
contain any properly embedded, essential, acylindrical, compact
orientable separating
surface. \end{Lemma}

A properly embedded, compact orientable surface
$(F^{2},\partial F^{2}) \hookrightarrow (M^{3}, \partial M^{3})$ is
acylindrical if
any embedded incompressible annulus
$(A^{2},\partial A^{2}) \hookrightarrow (M^{3}, F^{2})$, such that
$A^{2} \cap F^{2} = \partial A^{2}$, can be homotoped into $F^{2}$, by
a
homotopy supported in the side of the surface containing $A^{2}$.

By the annulus Theorem (cf. [Ja, Chap.VIII], [Sco]), a properly
embedded, essential
and acylindrical, compact, orientable and separating surface in
$M^{3}$
induces a malnormal (i.e. 1-acylindrical) amalgamated splitting of
$\pi_{1}(M^{3})$.
By a result of A. Karrass and D. Solitar [KS], a non-free,
two generated
group cannot admit such a malnormal splitting. \hfill$\Box$

To finish the proof of Corollary 3 when $X^{3}$ is hyperbolic, it
remains to show
that the properly embedded, essential, four punctured sphere
$(\Sigma^{2}, \partial \Sigma^{2}) \hookrightarrow (X^{3}, \partial
E(k))$ is
acylindrical in $X^{3}$.

If $(A^{2},\partial A^{2}) \hookrightarrow (X^{3}, \Sigma^{2})$, such
that $A^{2}
\cap \Sigma^{2} = \partial A^{2}$, is an incompressible annulus, then
the two components of
$\partial A^{2}$ are parallel essential simple curves on the four
punctured
sphere $\Sigma^{2}$. Otherwise, one component of $\partial A^{2}$ must
be boundary parallel on $\Sigma^{2}$. Since $k$ is a prime knot, the
other component of
$\partial A^{2}$ cannot be parallel to the boundary on $\Sigma^{2}$.
Hence, this curve must be bound a disk with two holes, that does not
contain the
other component of $\partial A^{2}$. This disk with two holes,
together with the
annulus $A^{2}$ gives an embedded three-punctured sphere, with
meridional boundary
components, in $E(k)$, which is impossible.

Now, the annulus $A^{2}$, together with the annulus on $\Sigma^{2}$
bounded by $\partial A^{2}$, give a torus $T^{2}$. This torus must be
compressible in the
side of $\Sigma^{2}$ that contains $A^{2}$, otherwise it would be
incompressible
in $X^{3}$, since $\Sigma^{2}$ is essential in $X^{3}$. Since $X^{3}$
is irreducible
and the components of $\partial A^{2}$ are not contained in a ball in
$X^{3}$, $T^{2}$ bounds a solid torus in one side of $\Sigma^{2}$,
which shows that the
annulus $A^{2}$ can be pushed on
$\Sigma^{2}$.
\hfill$\Box$


\footnotesize\begin{thebibliography}{XXXXX}

\bibitem[Ad]{Ad} C. Adams, \emph{Hyperbolic 3-manifolds with two
generators},
Comm. Anal. Geom. \textbf{4}, 1996, 181-206 181-206.

\bibitem[BBJ]{BBJ} D. Bachman, S.Bleiler and A. Jones, \emph{Two generator
3-manifolds},
in preparation.

\bibitem[BH]{BH} J.S.Birman and H. Hilden, \emph{Heegaard splittings
and branched coverings os
$S^3$}, Trans. Amer. Math. Soc. \textbf{213}, 1975, 315-352.

\bibitem[Ble1]{Ble1} S.Bleiler , \emph{Knots prime on many strings},
Trans. Amer. Math. Soc.
\textbf{282}, 1984, 385-401.

\bibitem[Ble2]{Ble2} S.Bleiler , \emph{Two generator cable knots are
tunnel
one},
Proc. Amer. Math. Soc.
\textbf{122}, 1994, 1285-1287.

\bibitem[Ble3]{Ble3} S.Bleiler , \emph{Little big knots},
Knot theory and its applications. Chaos, Solitons and Fractals \textbf{9},
No 4/5, 1998,
681-692.

\bibitem[BJ1]{BJ1} S.Bleiler and A. Jones, \emph{The free product of
groups with amalgamated
subgroup malnormal in a single factor}, J.Pure and Appl. Algebra
\textbf{127}, 1998, 119-136.

\bibitem[BJ2]{BJ2} S.Bleiler and A. Jones, \emph{On Two-Generator
satellite knots}, to appear in the Proceedings of the Stalling's conference.

\bibitem[BZg]{BZg} M.Boileau and H.Zieschang, \emph{Heegard genus of
closed orientable
Seifert 3-manifolds}, Invent. Math. \textbf{76}, 1984, 455-468.

\bibitem[BZn]{BZn} M.Boileau and B. Zimmermann, \emph{The
$\pi$-orbifold group of
a link}, Math.Z. \textbf{200}, 1989, 187-208.

\bibitem[BP]{BP} M. Boileau and J. Porti \emph{Geometrization of
3-orbifolds of cyclic type}, Ast\'eris\-que Monograph, \textbf{272}, 2001.

\bibitem[BS]{BS} F. Bonahon and L. Siebenmann \emph{The
characterisitc splitting
of irreducible compact 3-orbifolds}, Math. Math. \textbf{278}, 1987.

\bibitem[CHK]{CHK} D. Cooper, C. Hodgson and S. Kerchkoff
\emph{ Three-dimensional orbifolds and cone-manifolds},
 MSJ Memoirs \textbf{5}, 2000.

\bibitem[DD]{DD} W. Dicks and M. Dunwoody \emph{Groups acting on
graphs}, Cambridge
University Press, 1989.

\bibitem[Ep]{Ep} D.B.A. Epstein, \emph{Projective planes in
3-manifolds},
Proc. London Math. Soc. (3) \textbf{11}, 1961, 469-484.


\bibitem[EN]{EN} D.Eisenbud and W.Neumann {\em Three-dimensional link
theory and invariants of plane curve singularities},
 Ann. Math. Studies  \textbf{110}, Princeton University Press, 1985.


\bibitem[FRW]{FRW} B. Fine, G.Rosenberger and R.Weidmann,
\emph{Two-generated subgroups
of free products with commuting subgroups}, to appear in J. Pure Appl.
Algebra.

\bibitem[GL]{GL} C.McA. Gordon, R. Litherland
  {\em Incompressible planar surfaces in 3-manifolds}, Topology Appl.
\textbf{18}, 1984, 121-144.


\bibitem[GLM]{GLM} C.McA. Gordon, R. Litherland and K. Murasugi
  {\em Signatures of covering links}, Can. J. Math.
\textbf{33}, 1981, 331-394.

\bibitem[Gor]{Gor} C.McA. Gordon
  {\em Dehn surgery and satelitte knots}, Trans. Amer. Math. Soc.
\textbf{275}, 1983, 687-708.


\bibitem[Ha]{Ha} R. Hartley {\em Knots with free period}, Can. J.
Math.
\textbf{33}, 1981, 91-102.


\bibitem[He]{He} J.Hempel {\em 3-manifolds},
Ann. Math. Stud. \textbf{86}, 1976.

\bibitem[Ku]{Ku} M. Kuhn {\em Tunnelsof 2-bridge links.}  J. Knot
Theory Ramifications  \textbf{5}, 1996, 167-171.



\bibitem[Ja]{Ja} W.H.Jaco
{\em Lectures on three-manifold topology}, AMS, 1977.

\bibitem[JS1]{JS} W.H.Jaco and P.B.Shalen {\em Peripheral structure
of 3-manifolds}, Invent. Math. \textbf{38}, 1976, 55-87.

\bibitem[JS2]{JS} W.H.Jaco and P.B.Shalen {\em Seifert fibered spaces
in 3-manifolds},
Memoirs Amer. Math. Soc. \textbf{21} No. 220, 1979.

\bibitem[KW]{KW} I.Kapovich and R.Weidmann {\em Two-generated groups
acting on trees},
Arch. Math. \textbf{73}, 1999, 172-181.

\bibitem[KS]{KS} A. Karrass and D. Solitar {\em The free product of
two groups with a
malnormal amalgamated subgroup.} Can. J. Math. \textbf{23}, 1971,
933-959.

\bibitem[Ko1]{Ko1} T.Kobayashi {\em Non-separating incompressible
tori in 3-manifolds},
J. Math. Soc. Japan \textbf{36}, 1984, 11-22.

\bibitem[Ko2]{Ko2} T.Kobayashi {\em Structures of the {H}aken
manifolds with {H}eegaard
splittings of genus two}, Osaka J. Math. \textbf{21}, 1984, 437-455.

\bibitem[Ko3]{Ko3} T.Kobayashi {\em Structures of full {H}aken
manifolds},
Osaka J. Math. \textbf{24}, 1987, 173-215.

\bibitem[MP]{MP} J. McCool and A. Pietrowski {\em On free products
with amalgamation of
two infinite cyclic groups}, J. Algebra, \textbf{18}, 1971, 377-383.

\bibitem[Mon1]{Mon1} J.M. Montesinos {\em Variedades de Seifert que son
recubridadores
ciclicos
ramificados de dos hojas}
Bol. Soc. Mat. Mexicana \textbf{18}, 1973,
1-32.


\bibitem[Mon2]{Mon2} J.M. Montesinos {\em Classical tessellations and
three-manifolds},
 Springer-Verlag, Heidelberg, 1987.


\bibitem[MB]{MB} J. Morgan and H. Bass {\em The Smith conjecture},
Academic Press,
Orlando, Fl, 1984.

\bibitem[MSc]{MSc} Y.Moriah and J.Schultens {\em Heegaard splittings
of Seifert fibered spaces are either vertical or horizontal}, Topology
\textbf{37}, 1998,
1089-1112.

\bibitem[Mor]{Mor} K.Morimoto {\em On On minimum genus Heegaard
splittings of some orientable closed 3-manifolds},
Tokyo J. Math. \textbf{12}, 1989,
321-355.




\bibitem[MSa]{MSa} K.Morimoto and M.Sakuma {\em On unknotting tunnels
for knots}, Math. Ann. \textbf{289}, 1991,
143-167.



\bibitem[N]{N} F.H.Norwood {\em Every two-generator knot is prime}
Proc. Am. Math. Soc.
\textbf{86}, 1982, 143-147.

\bibitem[O]{O} P.Orlik {\em Seifert Manifolds}, Springer, Lect. Notes
in Math. \textbf{291}, 1972.


\bibitem[PRZ]{PRZ}
N.~Peczynski, G.~Rosenberger, and H.~Zieschang. {\em{\"U}ber
{E}rzeugende ebener
diskontinuierlicher {G}ruppen}, Invent. Math., \textbf{29}, 1975,
161--180.

\bibitem[Rei]{Rei}
A.~Reid
{\emph Some remarks on 2-generator hyperbolic 3-manifolds}, Lond.
Math. Soc. Lect. Note Ser. \textbf{173}, 1992, 209-219.

\bibitem[RiS]{RiS1} Y.~Rieck and E.~Sedgwick, \emph{Persistence of Heegaard
structures under
Dehn filling},
Topology and its Applications. \textbf{109}, 2001, 41-53.

\bibitem[RiS2]{RiS2} Y.~Rieck and E.~Sedgwick, \emph{Finiteness results for
Heegaard surfaces in
surgered manifolds},
Comm. Anal. Geom. \textbf{9}, 2001, 351-367.

\bibitem[RuS1]{RuS1}
H.~Rubinstein and M.~Scharlemann. {\em Comparing Heegaard
splittings-the
bounded case},
Trans. Amer. Math. Soc., \textbf{350}, 1998, 689--715.


\bibitem[RuS2]{RuS2}
H.~Rubinstein and M.~Scharlemann. {\em Genus two Heegaard splittings
of
orientable
three-manifolds} Geom. and Top. Monographs, \textbf{2}: Proceedings
of the
Kirbyfest, 1999,
489--553.


\bibitem[Sa1]{Sa1} M. Sakuma {\em Surface bundles over $S^1$ which
are 2-fold
branched cyclic coverings of $\mathbb S^3$.} Math. Sem. Notes, Kobe
Univ.
  \textbf{9}, 1981, 159-180.


\bibitem[Sa2]{Sa2} M. Sakuma {\em The Geometries of Spherical
Montesinos Links.} Kobe J.
Math. \textbf{7}, 1990, 167-190.

\bibitem[Scha]{Scha} M. Scharlemann {\em Tunnel number one knots
satisfy the Poenaru
conjecture} Topology and its Applications. \textbf{18}, 1984,
235-258.

\bibitem[Sc]{Sc} H. Schubert {\em Knoten mit zwei Br\"ucken} Math.Z.
\textbf{65},
1956, 133-170.

\bibitem[Sco]{Sco} P.Scott
{\em A new proof of the annulus and torus theorem}, A. J. Math.
\textbf{102}, 1980, 241-277.

\bibitem[Se]{Se} Z.Sela {\em Acylindrical accessibility for groups},
Invent. math.
\textbf{129}, 1997, 527-565.

\bibitem[S]{S} J.P.Serre {\em Trees}, New York, 1980

\bibitem[TO]{TO} M. Takahashi and M. Ochiai, {\em Heegaard diagrams
of torus bundles
over $S^1$}, Comm. Math. Univ. San. Pauli. \textbf{31}, 1982, 63-69.

\bibitem[Th]{Th} W. Thurston, {\em The geometry and topology of
3-manifolds} lecture notes, 1978.

\bibitem[Wa1]{Wa1} F.Waldhausen {\em Eine Klasse von 3-dimensionalen
Mannigfaltigkeiten I}
Invent. Math. \textbf{3}, 1967, 308-333.

\bibitem[Wa2]{Wa2} F.Waldhausen {\em Eine Klasse von 3-dimensionalen
Mannigfaltigkeiten II}
Invent. Math. \textbf{4}, 1967, 87-117.

\bibitem[W1]{W1} R.Weidmann {\em On the rank of amalgamated products
and product
knot groups}, Math. Ann. \textbf{312}, 1998, 761-771.

\bibitem[W2]{W2} R.Weidmann {\em Some 3-manifolds with 2-generated
fundamental group}, preprint

\bibitem[Z]{Z} H.Zieschang {\em {\"U}ber die Nielsensche
K\"urzungsmethode in freien
Produkten mit Amalgam}, Invent. Math. \textbf{10}, 1970, 4-37.

\bibitem[ZVC]{ZVC} H.Zieschang, E.Vogt and H.-D.Coldewey {\em
Surfaces and Planar
Discontinous Groups}, Springer, Lect. Notes in Math. \textbf{835},
1980.

\end{thebibliography}
\end{document}